\documentclass[a4paper,10pt,reqno]{amsart}
\usepackage[centertags]{amsmath}
\usepackage{amsfonts,amssymb,amsthm}
\usepackage[dvips]{graphicx}
\usepackage{psfrag}
\usepackage[english]{babel}
\usepackage{newlfont}
\usepackage{fancyhdr}
\newtheorem{theorem}{Theorem}[section]

\newtheorem{lemma}[theorem]{Lemma}

\newtheorem{corollary}[theorem]{Corollary}

\theoremstyle{definition}
\newtheorem{remark}[theorem]{Remark}
\newtheorem{definition}[theorem]{Definition}
\newtheorem{example}[theorem]{Example}

\numberwithin{equation}{section}

\newcommand{\Nb}{\mathbb{N}}
\newcommand{\NP}{\mathbb{N}^+}
\newcommand{\Rb}{\mathbb{R}}
\newcommand{\Op}{\mathcal{O}}
\newcommand{\E}{\mathcal{E}}
\newcommand{\F}{\mathcal{F}}

\newcommand{\LL}{\mathcal{L}}
\newcommand{\X}{\mathcal{X}}
\newcommand{\Y}{\mathcal{Y}}
\newcommand{\PP}{\mathcal{P}}
\newcommand{\UU}{\mathcal{U}}
\newcommand{\e}{\varepsilon}
\renewcommand{\O}{\Omega}
\renewcommand{\o}{\omega}
\newcommand{\G}{\Gamma}
\renewcommand{\a}{\alpha}

\newcommand{\cell}{\text{cell}}
\renewcommand{\hom}{\text{hom}}
\renewcommand{\div}{\text{div}}
\newcommand{\per}{\text{per}}
\newcommand{\Gt}{{\mathcal{G}_{\theta}}}
\newcommand{\Gtx}{{\mathcal{G}_{\theta(x)}}}
\newcommand{\dH}{\mathbf{d}_{\mathcal{H}}}
\newcommand{\dig}{\mathrm{diag}}

\renewcommand{\neg}{\negthinspace}

\newcommand{\ds}{\displaystyle}
\newcommand{\med}{- \hskip -1em \int}
\newcommand{\snorm}[1]{\left|#1\right|}
\newcommand{\la}{\langle}
\newcommand{\ra}{\rangle}
\makeatletter
\def\rightharpoonupfill@{\arrowfill@\relbar\relbar\rightharpoonup}
\newcommand{\xrightharpoonup}[2][]{\ext@arrow
0359\rightharpoonupfill@{#1}{#2}} \makeatother

\newcommand{\res}{\mathop{\hbox{\vrule height 7pt width .5pt depth 0pt
\vrule height .5pt width 6pt depth 0pt}}\nolimits}

\begin{document}
\title[A variational approach to the local character of G-closure]
{A variational approach to the local character of G-closure: the
convex case}
\author[J.-F. Babadjian \& M. Barchiesi]
{Jean-Fran\c cois Babadjian \& Marco Barchiesi}
\date{August 21, 2007}
\address{SISSA, Via Beirut 2-4, 34014 Trieste, Italy}
\email{babadjia@sissa.it} \email{barchies@sissa.it}
\maketitle

\begin{center}
\begin{minipage}{12cm}
\small{ \noindent {\bf Abstract.} This article is devoted to
characterize all possible effective behaviors of composite materials
by means of periodic homogenization. This is known as a $G$-closure
problem. Under convexity and $p$-growth conditions ($p>1$), it is
proved that all such possible effective energy densities obtained by
a $\Gamma$-convergence analysis, can be locally recovered by the
pointwise limit of a sequence of periodic homogenized energy
densities with prescribed volume fractions. A weaker locality result
is also provided without any kind of convexity assumption and the
zero level set of effective energy densities is characterized in
terms of Young measures. A similar result is given for cell
integrands which enables to propose new counter-examples to the
validity of the cell formula in the nonconvex case and to the
continuity of the determinant with respect
to the two-scale convergence.}\\

\small{ \noindent {\bf R\'esum\'e.} Cet article est consacr\'e \`a
la caract\'erisation de toutes les limites effectives possibles de
mat\'eriaux composites en terme d'homog\'en\'eisation p\'eriodique.
Ce probl\`eme est connu sous le nom de $G$-fermeture. Il est
d\'emontr\'e, sous des hypoth\`eses de convexit\'e et de croissance
$p>1$, que toutes ces densit\'es d'\'energies effectives, obtenues
lors d'une analyse par $\Gamma$-convergence, peuvent \^etre
localement vues comme la limite ponctuelle d'une suite de densit\'es
d'\'energies p\'eriodiquement homog\'en\'eis\'ees avec une fraction
de volume fix\'ee. Un r\'esultat plus faible est obtenu sans aucune
hypoth\`ese de convexit\'e et l'ensemble des z\'eros de ce type de
densit\'es d'\'energies effectives est caract\'eris\'e en terme de
mesures d'Young. Un r\'esultat similaire est donn\'e pour des
int\'egrandes cellulaires, permettant ainsi de proposer de nouveaux
contre-exemples quant \`a la validit\'e de la formule de cellule
dans le cas non convexe, ainsi qu'\`a la continuit\'e du
d\'eterminant par rapport \`a la convergence \`a double \'echelle.}

\vspace{10pt} \noindent {\bf Keywords:} $G$-closure, homogenization,
$\Gamma$-convergence, convexity, quasiconvexity, polyconvexity,
Young measures, two-scale convergence.

\vspace{6pt} \noindent {\bf 2000 M.S.C.:} 35B27, 35B40, 49J45,
73B27, 74E30, 74Q05.

\end{minipage}
\end{center}

\bigskip

\tableofcontents

\section{Introduction}

\noindent Composites are structures constituted by two or more
materials which are finely mixed at microscopic length scales.
Despite the high complexity of their microstructure, composites
appear essentially as homogeneous at macroscopic length scale. It
suggests to give a description of their effective properties as a
kind of average made on the respective properties of the
constituents. The \emph{Homogenization Theory} renders possible to
define properly such average, by thinking of a composite as a limit
(in a certain sense) of a sequence of structures whose
heterogeneities become finer and finer. There is a wide literature
on the subject; we refer the reader to \cite{Mil02} for a starting
point.

Many notions of convergence have been introduced to give a precise
sense to such asymptotic analysis. One of the most general is the
$H$-convergence (see \cite{MT,T4}), which permits to describe the
asymptotic behavior of a sequence of second order elliptic operators
in divergence form. This notion appears as a generalization of the
$G$-convergence (see \cite[Chapter 22]{DM}) introduced independently
for symmetric operators. Specifically, given a bounded open set
$\O \subseteq \Rb^n$, a sequence $\{A_k\}\subset L^\infty(\O;\Rb^{n
\times n})$ of uniformly elliptic tensors $H$-converges to
$A_\infty\in L^\infty(\O;\Rb^{n \times n})$ if, for every $g\in
H^{-1}(\O)$, the sequence of solutions $u_k$ of
$$-\div(A_k\nabla u)=g \text{ in } \mathcal D'(\O), \quad u \in H_0^1(\O)$$ satisfies
\begin{equation*}
\left\{\begin{array}{l} \displaystyle u_k\rightharpoonup u_\infty
\text{ weakly in } H_0^1(\O),
\\[5pt]
A_k\nabla u_k \rightharpoonup A_\infty\nabla u_\infty  \text{ weakly
in }L^2(\O;\Rb^n).
\end{array}\right.
\end{equation*}
\noindent \vspace{3pt} In particular, $u_\infty$ is the solution of
$$-\div(A_\infty \nabla u)=g \text{ in } \mathcal D'(\O), \quad u \in H_0^1(\O).$$
If the tensors $A_k$ describe a certain property (for instance the
conductivity) of the structures approximating the composite, it is
assumed that $A_\infty$ describes the effective behavior with
respect to that property.

Composite materials are characterized by three main features: the
different constituents (or phases), their volume fraction and their
geometric arrangement. A natural question which arises is the
following: given two constituents as well as their volume fraction,
what are all possible effective behaviors with respect to a certain
property? In the language of $H$-convergence this problem reads as
follows: given two tensors $A^{(1)}$ and $A^{(2)} \in \Rb^{n \times
n}$ corresponding to the conductivities of both constituents of the
composite and given a volume fraction $\theta \in [0,1]$, determine
all the tensors $A^* \in L^\infty(\O;\Rb^{n \times n})$ such that
there exists a sequence of characteristic functions $\{\chi_k\}
\subset L^\infty(\O;\{0,1\})$ (the geometry of the finer and finer
mixture) satisfying
\begin{equation*}\left\{\begin{array}{l}
A_{\chi_k} \text{ $H$-converges to }A^*,\\[0.2cm]
\displaystyle\med_\O\chi_k \, dx = \theta \;\; \text{ for all }k \in \Nb,
\end{array}\right.
\end{equation*}
where $A_{\chi_k}(x):=\chi_k(x) A^{(1)} + (1-\chi_k(x))A^{(2)}$.
The family of these effective tensors is called the $G$-closure of $\{A^{(1)},A^{(2)}\}$
with fixed volume fraction $\theta$.
Its determination is known as a {\it $G$-closure problem}
and is, in general, very difficult to solve except in some particular
cases that we shall briefly discuss later. Indeed, there is, in
general, no ``explicit'' formula to compute a $H$-limit. An
exception occurs in the periodic case, {\it i.e.}, when
\begin{equation*}
A_k(x)=A(\la k \,x \ra) \text{ for } x\in\O,
\end{equation*}
where $\la \cdot \ra$ denotes the fractional part of a vector
componentwise. In this very particular case it is possible to prove
(see \cite[Theorem 1.3.18]{A}) that the sequence $\{A_k\}$
$H$-converges to a constant tensor $A_\cell$ given by
\begin{equation}\label{Acell}
(A_\cell)_{ij}=\int_Q A(y)[e_i+\nabla\phi_i(y)]\cdot
[e_j+\nabla\phi_j(y)]\, dy \quad \text{ for }\quad i,j\in\{1,\ldots,
n\},
\end{equation}
where $\{e_1,\ldots,e_n\}$ is the canonical basis of $\Rb^n$,
$Q:=(0,1)^n$ and $\phi_i$ is a solution (unique up to an additive
constant) of the problem $$-\div \big( A[e_i+\nabla\phi]\big) =0
\text{ in }\mathcal D'(Q),\quad \phi \in H^1_\per(Q).$$

With the aim of finding optimal bounds and an analytical description of the
$G$-closure set of two isotropic conductors in dimension two, it has
been proved independently in \cite{LC} and \cite{T1} a locality
property which underlines the importance of periodic structures. It
states that every effective tensors, obtained by mixing two
materials with a volume fraction $\theta$, can be locally
recovered by the pointwise limit of a sequence of effective tensors,
each of them obtained by a periodic mixture with the same proportion $\theta$.
In other words, in that case, periodic
mixtures capture every kind of mixtures, which enables to reduce the
study of such problems to periodic geometries and materials with
homogeneous effective behavior, since periodic homogenization
produces homogeneous limiting ones. This locality result has been
subsequently generalized by {\sc Dal Maso} and {\sc Kohn} in an
unpublished work to higher dimension and not necessarily isotropic
conductors (a proof can be found {\it e.g.} in \cite[Chapter~2]{A}).
Given two tensors $A^{(1)}$ and $A^{(2)} \in \Rb^{n \times n}$,
for any $\theta \in [0,1]$, define
\begin{eqnarray*}
P_\theta & := & \Big\{ A^* \in \Rb^{n \times n}: \text{ there
exists } \chi \in L^\infty(Q;\{0,1\})\\[-0.2cm]&& \quad\text{ such that }
\int_Q \chi\, dx =\theta \text{ and } A^*= (A_\chi)_\cell\Big\},
\end{eqnarray*}
where $A_\chi(x):=\chi(x) A^{(1)} + (1-\chi(x))A^{(2)}$ and
$(A_\chi)_\cell$ is defined as in \eqref{Acell}. Denote by
$G_\theta$ the closure of $P_\theta$ in $\Rb^{n \times n}$. The
local representation of $G$-closure states that if $A^* \in
L^\infty(\O;\Rb^{n \times n})$ and $\theta \in L^\infty(\O;[0,1])$,
then $A^*(x_0) \in G_{\theta(x_0)}$ for a.e. $x_0 \in \O$ if and
only if there exists a sequence of characteristic functions
$\{\chi_k\} \subset L^\infty(\O;\{0,1\})$ such that
\begin{equation*}\left\{\begin{array}{l}
A_{\chi_k} \text{ $H$-converges to }A^*,\\[0.2cm]
\chi_k \rightharpoonup \theta \text{ weakly* in }L^\infty(\O;[0,1]).
\end{array}\right.
\end{equation*}
Later, a further generalization has been given in \cite{R} to the
case of nonlinear elliptic and strictly monotone operators in
divergence form, using an extension of the $H$-convergence to the
nonlinear monotone setting which can be found in \cite{CDMD}.

\vspace{5pt}

In this paper, we wish to study similar $G$-closure problems in the
framework of nonlinear elasticity. Of course, this degree of
generality makes it quite impossible to characterize analytically a
$G$-closure set. However, it is still possible to formulate a
locality result which finds a similar statement also in that case.
As nonlinear elasticity rests on the study of equilibrium states, or
minimizers, of a suitable energy, we fall within the framework of
the Calculus of Variations and it is natural to use the notion of
$\G$-convergence (see \cite{B06}) to describe the effective
properties of a composite. It is a variational convergence which
permits to analyze the asymptotic behaviour of a sequence of
minimization problems of the form
\begin{equation}\label{one}
\mathrm{m}_k=\mathrm{min}\lbrace F_k(u) \,:\, u\in\UU\rbrace,
\end{equation}

\vspace{4pt} \noindent where $\UU$ is a suitable functional space
(in the sequel $W^{1,p}(\O;\Rb^m)$ with $p\in(1,+\infty)$) and
$F_k$ is a sequence of functionals representing in our case the
energies of the structures approximating the composite. The
$\G$-limit $F_\infty$ of the sequence $F_k$ has the following
fundamental property: if $u_k$ is a solution of the minimization
problem \eqref{one} and $u_k\rightarrow u_\infty$ in $\UU$, then
$F_k(u_k)\rightarrow F_\infty(u_\infty)$ and $u_\infty$ is a
solution of the minimization problem
\begin{equation*}
\mathrm{m}_\infty=\mathrm{min}\lbrace F_\infty(u) \,:\, u\in\UU\rbrace.
\end{equation*}
This fact supports the assumption that $F_\infty$ describes in a good way the effective
energy of the composite.

\vspace{5pt}

Let us briefly explain our study in the special case of mixtures of
two different materials. Let $W^{(1)}$ and $W^{(2)}:\Rb^{m \times n}
\to [0,+\infty)$ be two functions satisfying standard $p$-growth and
$p$-coercivity conditions. The functions
$W^{(1)}$ and $W^{(2)}$ stand for the stored energy densities of two
homogeneous nonlinearly elastic materials. If $\chi\subset
L^\infty(\O;\{0,1\})$ is a characteristic function, we use the
notation
$W_\chi(x,\xi):=\chi(x)W^{(1)}(\xi)+(1-\chi(x))W^{(2)}(\xi)$. The
function $W_\chi$ can be thought of as the stored energy density of
a composite material obtained as a mixture of both previous ones.

Without any kind of convexity assumption on $W^{(1)}$ and $W^{(2)}$,
we prove a weaker locality property about effective energy densities
(Theorem \ref{localiz}) which states that if $W^*:\O \times \Rb^{m
\times n} \to [0,+\infty)$ is a Carath\'eodory function and $\theta
\in L^\infty(\O;[0,1])$, then the following conditions are
equivalent:
\begin{itemize}
\item[i)]
there exists a sequence of characteristic functions
$\{\chi_k\} \subset L^\infty(\O;\{0,1\})$ such that
\begin{equation}\label{dossier}
\left\{\begin{array}{l} \displaystyle
\int_\Omega W^*(x,\nabla u(x))\, dx = \Gamma\text{-}\lim_{k \to
+\infty} \int_\O
W_{\chi_k}(x,\nabla u(x))\, dx,\\[0.4cm]
\chi_k \rightharpoonup \theta \text{ weakly* in }L^\infty(\O;[0,1]);
\end{array}\right.
\end{equation}
\item[ii)]
for a.e. $x_0 \in \O$, there exists a sequence of characteristic
functions $\{\tilde \chi_k\} \subset L^\infty(Q;\{0,1\})$ such that
\begin{equation*}
\left\{\begin{array}{l} \displaystyle \int_Q W^*(x_0,\nabla u(y))\, dy
= \Gamma\text{-}\lim_{k \to +\infty} \int_Q
W_{\tilde \chi_k}(y,\nabla u(y))\, dy,\\[0.4cm]
\tilde \chi_k \rightharpoonup \theta(x_0) \text{ weakly* in }L^\infty(Q;[0,1]).
\end{array}\right.
\end{equation*}
\end{itemize}
We interpret this result as being enough to consider homogeneous
effective energies with asymptotically homogeneous microstructure.
In Theorem \ref{ponte} we perform another characterization of the
effective energy densities through a fine analysis of their zero
level set. It is realized by means of gradient Young measures
supported in the union of the zero level sets of $W^{(1)}$ and
$W^{(2)}$.

When the functions $W^{(1)}$ and $W^{(2)}$ are both convex, we can carry out
a deeper locality result. As before, for any $\theta \in [0,1]$, we
define the set $P_\theta$ made of all effective energy densities
obtained by a periodic homogenization of a mixture of $W^{(1)}$ and
$W^{(2)}$ in proportions $\theta$ and $1-\theta$:
\begin{eqnarray*}P_\theta & := & \Big\{W^*:\Rb^{m \times n} \to
[0,+\infty): \text{ there exists }\chi \in L^\infty(Q;\{0,1\})\\
&&\hspace{1cm}\ds\text{such that } \int_Q \chi\, dx =\theta \text{
and } W^*=(W_\chi)_\hom  \Big\},
\end{eqnarray*}
where
\begin{equation}\label{intro hom}
(W_\chi)_{\hom}(\xi):=\inf_{j\in\NP}\inf\,\Bigl\lbrace\med_{(0,j)^n}
W_\chi\bigl(\langle y\rangle, \xi+ \nabla\phi(y)\bigr)\, dy \,:\,
\phi\in W^{1,p}_{\per}\bigl((0,j)^n;\Rb^m\bigr)\Bigr\rbrace.
\end{equation}

\vspace{3pt} \noindent The function $(W_\chi)_{\hom}$ is called
\emph{homogenized integrand} associated to $W_\chi$ and it is
possible to prove (independently of the convex assumption) that if
$\chi_k(x)=\chi(\la k \,x \ra)$ for $x\in\O$, then
$$\int_\Omega (W_\chi)_{\hom}(\nabla u(x))\, dx = \Gamma\text{-}\lim_{k \to
+\infty} \int_\O W_{\chi_k}(x,\nabla u(x))\, dx.$$

In general, the set $P_\theta$ is not closed for the pointwise
convergence so that it is natural to consider its closure denoted as
before by $G_\theta$.
Indeed, there exist effective energies that are not exactly reached by a periodic
mixture: thanks to Theorem \ref{ponte}, we built an example in the case
of a mixture of four materials (Example \ref{tartar}).
Our main result, Theorem \ref{BB}, states that if $W^*:\O \times \Rb^{m \times n} \to
[0,+\infty)$ is a Carath\'eodory function and $\theta \in
L^\infty(\O;[0,1])$, then there exists a sequence of characteristic
functions $\{\chi_k\} \subset L^\infty(\O;\{0,1\})$ satisfying
\eqref{dossier} if and only if $W^*(x_0,\cdot) \in G_{\theta(x_0)}$
for a.e. $x_0 \in \O$.

This result was stated in \cite[Conjecture 3.15]{FF} and it was
useful to study a quasi-static evolution model for the interaction
between fracture and damage. Indeed the authors used this conjecture
to prove the wellposedness of an incremental problem at fixed time
step.

If we assume the functions $W^{(1)}$ and $W^{(2)}$ to be
differentiable, one can write the Euler-Lagrange equation associated
to the minimization problem. It is a nonlinear elliptic partial
differential equation in divergence form and the first order
convexity condition means that the operators $$\xi \mapsto
\frac{\partial W^{(1)}}{\partial \xi}(\xi)\quad \text{ and }\quad
\xi \mapsto \frac{\partial W^{(2)}}{\partial \xi}(\xi)$$ are
monotone. Hence our locality result reduces to a particular case of
\cite{R}. Furthermore, if $W^{(1)}(\xi)=A^{(1)}\xi \cdot \xi$ and
$W^{(2)}(\xi)=A^{(2)}\xi \cdot \xi$ for some symmetric matrices
$A^{(1)}$ and $A^{(2)} \in \Rb^{n \times n}$, then Theorem \ref{BB}
is nothing but {\sc Dal Maso} and {\sc Kohn}'s result in the special
case of symmetric operators. Note that since $H$-convergence
problems have, in general, no variational structure, our result do
not generalize the preceding ones, but rather gives a variational
version of them.

Our proof of Theorem \ref{BB} works only when the functions
$W^{(1)}$ and $W^{(2)}$ are both convex. The nonconvex case seems to
be much more delicate to address. The reason is that the
homogenization formula \eqref{intro hom} in the convex case is reduced
to a single cell formula, \textit{i.e.}, $(W_\chi)_{\hom}=(W_\chi)_{\cell}$,
where
\begin{equation*}
(W_\chi)_{\cell}(\xi):=\inf\,\Bigl\lbrace\int_Q W_\chi\bigl(y, \xi+
\nabla\phi(y)\bigr)\, dy \,:\, \phi\in
W^{1,p}_{\per}\bigl(Q;\Rb^m\bigr)\Bigr\rbrace.
\end{equation*}
On the contrary, it is known that both formula do not coincide even
in the quasiconvex case because of the counter-example \cite[Theorem~4.3]{M}.
To illustrate this phenomena, we present here another much simpler
counter-example, Example \ref{unico}, based on a rank-one connection
argument and on the characterization of the zero level set of cell
integrands under quasiconvex assumption, Lemma \ref{zero energy}.
Indeed we prove that when $m=n=2$, there exist suitable densities
$W^{(1)}$ and $W^{(2)}$ (one of them being non convex) and a
suitable geometry $\chi$ such that
$$(W_\chi)_\cell(C)>0 \;\text{ while }\; (W_\chi)_\hom(C)=0$$
for some matrix $C\in \Rb^{2 \times 2}$. Moreover, exploiting the
lower semicontinuity property of certain integral functionals with
respect to the two-scale convergence and noticing that the functions
$W^{(1)}$ and $W^{(2)}$ in the previous example can be taken
polyconvex, we also show that, in general, one cannot expect the
determinant to be continuous with respect to the two-scale
convergence (Example \ref{no-unico}). \vspace{5pt}

The paper is organized as follows: in Section \ref{notations}, we
introduce notations used in the sequel and we recall basic facts
about $\G$-convergence, periodic homogenization, Young measures and
two-scale convergence. In Section \ref{localprop}, we prove a
locality property, Theorem \ref{localiz}, enjoyed by effective
energy densities obtained by mixing $N$ different materials.
Section~\ref{mesyoung} is devoted to a fine analysis of the zero
level set of such effective energy densities by means of Young
measures: this is achieved in Theorem \ref{ponte}. We also perform a
similar analysis for cell kind integrands. In Section \ref{gcl}, we
give a proof of our main result, Theorem \ref{BB}, which states the
local property with respect to the set $G_\theta$ in the case of
convex stored energy densities. Moreover, thanks to the previous
Young measure analysis, we show by an example based on the Tartar
square (Example \ref{tartar}) that in general $P_\theta \subsetneq
G_\theta$. Finally, we present in Section \ref{ce} a new
counter-example to the validity of the cell formula in the nonconvex
case (Example~\ref{unico}) and to the continuity of the determinant
with respect to the two-scale convergence (Example~\ref{no-unico}).

\vspace{10pt}


\section{Notations and preliminaries}\label{notations}

\noindent Throughout the paper, we employ the following notations:
\begin{itemize}
\item $\O$ is a bounded open subset of $\Rb^n$;
\item $Q$ is the unit cell $(0, 1)^n$;
\item if $a \in \Rb^n$ and $\rho>0$, $B_\rho(a)$ stands for
the open ball in $\Rb^n$ of center $a$ and radius $\rho$ while
$Q_\rho(a):=a+(0,\rho)^n$ is the open cube in $\Rb^n$ with corner
$a$ and edge length $\rho$;
\item $\Op(\O)$ denotes the family of all open subsets of $\O$;
\item $\LL^n$ denotes the $n$-dimensional Lebesgue measure;
\item the symbol $\med_E$ stands for the average
$\LL^n(E)^{-1}\int_E$;
\item if $\mu$ is a Borel measure in $\Rb^s$ and $E \subseteq \Rb^s$
is a Borel set, the measure $\mu\res\, E$ stands for the restriction
of $\mu$ to $E$, {\it i.e.} for every Borel set $B \subseteq \Rb^s$,
$(\mu\res\, E)(B)=\mu(E\cap B)$;
\item $W^{1,p}_\per(Q;\Rb^m)$ is the space of $W^{1,p}_\textrm{loc}(\Rb^n;\Rb^m)$ functions which are $Q$-periodic;
\item given $x\in \Rb^n$, $\la x \ra$ stands for the fractional part of $x$ componentwise;
\item $\rightharpoonup$ (resp.
$\xrightharpoonup[]{*}$) always denotes weak (resp. weak*) convergence.
\end{itemize}

\subsection{$\G$-convergence}

\noindent Let $\a\in\Rb^3$ with $\a_1,\a_2,\a_3>0$ and
$p\in(1,+\infty)$. We denote by $\F(\a,p)$ the set of all continuous
functions $f : \Rb^{m\times n} \to [0,+\infty)$ satisfying the
following coercivity and growth conditions:
\begin{equation}\label{pgrowth}
\a_1\snorm{\xi}^p-\a_2 \leq f(\xi) \leq \a_3 (1+\snorm{\xi}^p)
\quad\text{for any $\xi \in \Rb^{m \times n}$}.
\end{equation}
Moreover, we denote by $\F(\a,p,\O)$ the set of all Carath\'eodory
functions $f :\O\times\Rb^{m\times n} \to [0,+\infty)$ such that
$f(x,\cdot)\in \F(\a,p)$ for a.e. $x\in\O$.

We recall the definition of $\G$-convergence, referring to
\cite{BD,B06,DM} for a comprehensive treatment on the subject.

\begin{definition}\label{defgammaconv}
Let $f$, $f_k\in\F(\a,p,\O)$ and define the functionals $F$,
$F_k:W^{1,p}(\O;\Rb^m)\to [0,+\infty)$ by
\begin{equation*}
F(u):=\int_\O f(x,\nabla u)\, dx\quad \text{ and } \quad
F_k(u):=\int_\O f_k(x,\nabla u)\, dx.\end{equation*} We say that the
sequence $\{F_k\}$ $\G$-converges to $F$ (with respect to the weak
topology of $W^{1,p}(\O;\Rb^m)$) if for every $u\in W^{1,p}(\O;\Rb^m)$
the following conditions are satisfied:
\begin{itemize}
\item[i)]\emph{liminf inequality}: for every sequence $\{u_k\} \subset W^{1,p}(\O;\Rb^m)$
such that $u_k \rightharpoonup u$ in $W^{1,p}(\O;\Rb^m)$,
\begin{equation*}
\hspace*{-2pt}F(u)\leq\liminf_{k\rightarrow+\infty}F_k(u_k);
\end{equation*}
\item[ii)]\emph{recovery sequence}: there exists a sequence $\{u_k\} \subset W^{1,p}(\O;\Rb^m)$
such that $u_k\rightharpoonup u$ in $W^{1,p}(\O;\Rb^m)$ and
\begin{equation*}
F(u)=\lim_{k\rightarrow+\infty}F_k(u_k).
\end{equation*}
\end{itemize}
\end{definition}

In the following remark, we state without proof standard results on
$\G$-convergence of integral functionals which can be found in the
references given.
\begin{remark}
\begin{itemize}\label{G}
\item[]
\item[i)] Thanks to the coercivity condition \eqref{pgrowth},
Definition \ref{defgammaconv} coincides with the standard definition
given in a topological space (see \cite[Proposition 8.10]{DM}).
\item[ii)] If $\{f_k\}$ is a sequence in $\F(\a,p,\O)$, then there exist a subsequence $\{f_{k_j}\}$ and a
function $f\in\F(\a,p,\O)$ such that $\{F_{k_j}\}$ $\G$-converges to
$F$ (see \cite[Theorem 12.5]{BD}).
\item[iii)] If $\{F_k\}$ $\G$-converges to $F$, then $F$ is sequentially weakly lower
semicontinuous in $W^{1,p}(\O;\Rb^m)$ and therefore $f$ is
quasiconvex with respect to the second variable (see
\cite[Statement~II-5]{AF}). Moreover, for every $U\in\Op(\O)$ and
every $u \in W^{1,p}(U;\Rb^m)$, we have
\begin{equation*}
\int_U f(x,\nabla u)\, dx=\G\text{-}\lim_{k \to +\infty} \int_U
f_k(x,\nabla u)\, dx.
\end{equation*}
\item[iv)] If $\{F_k\}$ $\G$-converges to $F$, defining
$\widetilde{F}_k:W^{1,p}(\O;\Rb^m)\rightarrow[0,+\infty)$ by
\begin{equation*}
\widetilde{F}_k(u):=\int_\O Qf_k(x,\nabla u)\, dx,
\end{equation*}
then $\{\widetilde{F}_k\}$ $\G$-converges to $F$ as well, where
$Qf_k(x,\cdot)$ denotes the quasiconvex envelope of $f_k(x,\cdot)$.
\item[v)] Given $f$ and $f_k\in\F(\a,p,\O)$, if the functions $f_k(x,\cdot)$ are quasiconvex
for a.e. $x \in \O$ and for all $\xi \in \Rb^{m \times n}$
\begin{equation*}
f_k(\cdot,\xi) \to f(\cdot,\xi) \quad\text{pointwise a.e. in $\O$},
\end{equation*}
then $\{F_k\}$ $\G$-converges to $F$ (see \cite[Lemma~7 and Remark~15]{Gl});
\item[vi)] If $\{u_k\} \subset W^{1,p}(\O;\Rb^m)$ is a
sequence such that $u_k \rightharpoonup u$ in $W^{1,p}(\O;\Rb^m)$
and $F_k(u_k) \to F(u)$, then one can assume that $u_k=u$ in a
neighborhood of $\partial \O$ (see \cite[Proposition~11.7]{BD}).
\end{itemize}
\end{remark}
The following technical result will be useful in the analysis of the
local properties of mixtures in Section \ref{localprop}.
\begin{lemma}\label{rho}
Let $f\in\F(\a,p,\O)$ be quasiconvex with respect to the second
variable. There exists a set $Z \subseteq \O$ with $\LL^n(Z)=0$ such
that for every $\{\rho_j\} \searrow 0^+$ and every $x_0 \in  \O
\setminus Z$, $$\G\text{-}\lim_{j\to +\infty}\int_Q f(x_0 + \rho_j y
,\nabla u(y))\, dy=\int_Q f(x_0,\nabla u(y))\, dy \quad \text{ for
all }u \in W^{1,p}(Q;\Rb^m).$$
\end{lemma}
\begin{proof}
By \cite[Lemma~5.38]{AFP}, there exists a $\LL^n$-negligible set $Z
\subseteq \O$ such that for any $x_0 \in \O\setminus Z$ and any
sequence $\{\rho_j\} \searrow 0^+$, one can find a subsequence
$\{\rho_{j_k}\}$ of $\{\rho_j\}$ and a set $E\subseteq Q$ with
$\LL^n(E)=0$ such that $$f(x_0+\rho_{j_k} y,\xi) \to f(x_0,\xi)$$
for every $\xi \in \Rb^{m \times n}$ and every $y \in Q \setminus
E$. Hence from Remark \ref{G} (v), one gets that
$$\G\text{-}\lim_{k \to + \infty}\int_Q f(x_0+\rho_{j_k}y,\nabla u(y))\, dy =
\int_Qf(x_0,\nabla u(y))\, dy, \text{ for every } u \in
W^{1,p}(Q;\Rb^m).$$ Since the $\G$-limit does not depend upon the
choice of the subsequence, we conclude in light of \cite[Proposition
~8.3]{DM} that the whole sequence $\G$-converges.
\end{proof}

\subsection{Periodic homogenization}

\noindent We now recall standard results concerning periodic
homogenization. We call \emph{cell integrand} (resp.
\emph{homogenized integrand}) associated to $f\in\F(\a,p,Q)$, the
function $f_{\cell}\in\F(\a,p)$ (resp. $f_{\hom}\in\F(\a,p)$)
defined by
\begin{equation}\label{fcell}
f_{\cell}(\xi):=\inf \Bigl\lbrace\int_Q f(y, \xi+\nabla \phi)\, dy
\;:\; \phi\in W^{1,p}_{\per}(Q;\Rb^m)\Bigr\rbrace
\end{equation}
\begin{equation*}
\left(\text{resp.
}f_{\hom}(\xi):=\inf_{j\in\NP}\inf\,\Bigl\lbrace\med_{(0,j)^n}
f\bigl(\langle y\rangle, \xi+\nabla \phi\bigr)\, dy \,:\, \phi\in
W^{1,p}_{\per}\bigl((0,j)^n;\Rb^m\bigr)\Bigr\rbrace\right).
\end{equation*}

\begin{remark}\label{minfcell}
If $f$ is quasiconvex in the second variable, then from standard
lower semicontinuity results
$$f_{\cell}(\xi)=\min \Bigl\lbrace\int_Q f(y, \xi+\nabla \phi)\, dy
\;:\; \phi\in W^{1,p}_{\per}(Q;\Rb^m)\Bigr\rbrace.$$
\end{remark}

\vspace{6pt}

\noindent
The following theorem is well known in the theory of
$\G$-convergence (See \cite[Section 14]{BD}, \cite{B85} or \cite{M}).
\begin{theorem}\label{BM}
Let $f\in\F(\a,p,Q)$ and $\{\e_k\} \searrow 0^+$, then the
functionals $F_k:W^{1,p}(\O;\Rb^m)\rightarrow[0,+\infty)$ defined by
\begin{equation*}
F_k(u):=\int_\O f\left(\left\langle\frac{x}{\e_k}\right\rangle,
\nabla u\right)\, dx
\end{equation*}
$\G$-converges to $F_\hom:W^{1,p}(\O;\Rb^m) \to [0,+\infty)$, where
$$F_{\hom}(u):=\int_\O f_{\hom}(\nabla u)\, dx.$$ In particular
$f_{\hom}$ is a quasiconvex function. Under the additional
hypothesis that $f$ is convex with respect to the second variable,
then $f_{\hom}=f_{\cell}$.
\end{theorem}

\subsection{Young measures}

\noindent
Let $s \in \Nb$, we denote by $\PP(\Rb^s)$ the space
of probability measures in $\Rb^s$ and by $\Y(\O;\Rb^s)$ the space
of maps $\mu:x\in\O \mapsto \mu_x\in \PP(\Rb^s)$ such that the
function $x \mapsto \mu_x(B)$ is Lebesgue measurable for every Borel
set $B\subseteq\Rb^s$ (see \cite[Definition~2.25]{AFP}).

\vspace{6pt} The following result is a version of the Fundamental
Theorem on Young Measures (a proof can be found in
\cite{B88,Mul99bis}). Under suitable hypothesis, it shows that the
weak limit of a sequence of the type $\{f\left(\cdot,
w_k(\cdot)\right)\}$ can be expressed through a suitable map
$\mu\in\Y(\O;\Rb^s)$ associated to $\{w_k\}$.

\begin{theorem}\label{fund young}
Let $\{w_k\}$ be a bounded sequence in $L^1(\O;\Rb^s)$. There exist a
subsequence $\{w_{k_j}\}$ and a map $\mu\in\Y(\O;\Rb^s)$ such that
the following properties hold:
\begin{enumerate}
\item[i)] if $f:\O\times\Rb^s\rightarrow[0,+\infty)$ is a Carath\'eodory function, then
\begin{equation*}
\liminf_{j\to +\infty}\int_\O f\bigl(x, w_{k_j}(x)\bigr)\,
dx\geq\int_\O\overline{f}(x)\, dx
\end{equation*}
where
\begin{equation*}
\overline{f}(x):=\int_{\Rb^s}f(x,\xi)\, d\mu_x(\xi);
\end{equation*}
\item[ii)] if $f:\O\times\Rb^s\to \Rb$ is a Carath\'eodory function and
$\{f(\cdot, w_k(\cdot))\}$ is equi-integrable, then $f(x,\cdot)$ is
$\mu_x$-integrable for a.e. $x\in\O$, $\overline{f}\in L^1(\O)$ and
\begin{equation*}
\lim_{j \to +\infty}\int_\O f\bigl(x, w_{k_j}(x)\bigr)\,
dx=\int_\O\overline{f}(x)\, dx;
\end{equation*}
\item[iii)] if $\mathfrak{A}\subseteq\Rb^s$ is compact, then $\mathrm{supp}\,\mu_x\subseteq\mathfrak{A}$ for a.e. $x\in\O$ if and only if $\mathrm{dist}(w_{k_j},\mathfrak{A})\rightarrow 0$ in measure.
\end{enumerate}
\end{theorem}
\begin{definition}
The map $\mu\in\Y(\O;\Rb^s)$ is called the \emph{Young measure}
generated by $\{w_{k_j}\}$.
\end{definition}

\begin{remark}\label{1049}
We denote by $\overline{\mu}_x:=\int_{\Rb^s}\xi\,d\mu_x$ the
\emph{barycenter} of $\mu$. If the sequence $\{w_k\}$ is
equi-integrable and generates $\mu$, then $w_k \rightharpoonup
\overline{\mu}$ in $L^1(\O;\Rb^s)$. To see this, it is sufficient to
apply Theorem \ref{fund young} (ii) with $f(x,\xi)=\varphi(x)\xi^{(i)}$
where $\varphi\in L^\infty(\O)$ and $\xi^{(i)}$ is the $i$'th component
of $\xi$ ($i\in\{1,\ldots,s\}$).
\end{remark}

In the sequel we are interested in Young measures generated by
sequences of gradients. We refer the reader to \cite{Mul99bis,Pe97}
and references therein for a more exhaustive study.
\begin{definition}
A map $\nu\in\Y(\O;\Rb^{m\times n})$ is called \emph{gradient Young
measure} if there exists
a bounded sequence $\{u_k\}$ in $W^{1,p}(\O;\Rb^m)$ such that $\{\nabla u_k\}$ generates $\nu$.\\
A probability measure $\sigma\in\PP(\Rb^{m\times n})$ is called
\emph{homogeneous gradient Young measure} if there exists a gradient
Young measure $\nu$ such that $\sigma=\nu_x$ for a.e. $x\in\O$.
\end{definition}

\begin{remark}
\label{Y} If $\nu$ is a gradient Young measure, then for a.e.
$x\in\O$ the measure $\nu_x$ is a homogeneous gradient Young measure
(see \cite[Theorem 8.4]{Pe97}).
\end{remark}

\subsection{Two-scale convergence}

\noindent We here briefly gather the definition and some of the main
results about two-scale convergence.
\begin{definition}
Let $\{w_k\}$ be a bounded sequence in $L^p(\O;\Rb^s)$ and
$\{\e_k\} \searrow 0^+$. We say that $\{w_k\}$ two-scale converges
to a function $w=w(x, y)\in L^p(\O\times Q;\Rb^s)$ (with respect to
$\{\e_k\}$) if
\begin{equation*}
\lim_{k\rightarrow+\infty}\int_{\O}\varphi(x)\,
\phi\left(\frac{x}{\e_k}\right) w_k(x)\, dx =\int_{\O\times
Q}\varphi(x)\, \phi\left(y\right) w\left(x,y\right)\, dx\,dy
\end{equation*}
for any $\varphi\in C_c^\infty(\O)$ and any $\phi\in
C_\per^\infty(Q)$; we simply write $w_k\rightsquigarrow w$.
\end{definition}

\begin{remark}\label{2scale}
\begin{itemize}
\item[]
\item[i)] Every bounded sequence in $L^p(\O;\Rb^s)$ admits a
two-scale converging subsequence (see \cite[Theorem~7]{LNW}).
\item[ii)] Let $f:Q \times \Rb^s \to [0,+\infty)$ be a Carath\'eodory
function convex with respect to the second variable and $\{\e_k\}
\searrow 0^+$. If $w_k \rightsquigarrow w$ (with respect to
$\{\e_k\}$), then
\begin{equation*}
\liminf_{k \to +\infty} \int_\O f \left(\left\langle
\frac{x}{\e_k}\right\rangle,w_k(x) \right)dx \geq \int_{\O \times Q}
f(y,w(x,y))\, dx\, dy.
\end{equation*}
This result is a direct consequence of \cite[Theorem~4.14]{Bar05}
together with Jensen's Inequality.
\item[iii)] Let $\{u_k\}$ be a sequence weakly converging in $W^{1,p}(\O;\Rb^m)$ to
a function $u$. Then $u_k \rightsquigarrow u$ and there exist a
subsequence (not relabeled) and $v \in
L^p(\O;W^{1,p}_\per(Q;\Rb^m))$ such that $\nabla u_k
\rightsquigarrow \nabla u + \nabla_y v$ (see
\cite[Theorem~13]{LNW}).
\end{itemize}

\end{remark}


\section{Local properties of mixtures}\label{localprop}

\noindent Let $W^{(1)},\ldots,W^{(N)}$ belong to $\F(\a,p)$. We
interpret these functions as the stored energy densities of $N$
different nonlinearly elastic materials and we are interested in
mixtures between them.

\begin{definition}\label{characfunct}
Let $\X(\O)$ be the family of all functions
$\chi=(\chi^{(1)},\ldots,\chi^{(N)}) \in L^\infty(\O;\{0,1\}^N)$
such that $\sum_{i=1}^N\chi^{(i)}(x)=1$ in $\O$. Equivalently,
$\chi\in\X(\O)$ if there exists a measurable partition $\lbrace
P^{(i)}\rbrace_{i=1,\ldots,N}$ of $\O$ such that
$\chi^{(i)}=\chi_{P^{(i)}}$ for $i=1,\ldots,N$.
\end{definition}

If $\chi\in \X(\O)$, we define $W_\chi\in\F(\a,p,\O)$ by
$$W_\chi(x,\xi):=\sum_{i=1}^N \chi^{(i)}(x) W^{(i)}(\xi), \quad \text{for every }
x \in \O \text{ and }\xi \in \Rb^{m \times n}$$ and the functional
$F_\chi : W^{1,p}(\O;\Rb^m) \to [0,+\infty)$ by
\begin{equation*}
F_\chi(u):=\int_\O W_{\chi}(x,\nabla u)\, dx.
\end{equation*}

\noindent The function $W_\chi$ can be thought of as the stored
energy density of a composite material obtained by mixing $N$
components having energy density $W^{(i)}$ and occupying the
reference configuration $P^{(i)}$ in $\O$ at rest. The elastic
energy of this composite material under a deformation $u \in
W^{1,p}(\O;\Rb^m)$ is given by $F_\chi(u)$.

To describe efficiently the complexity of a composite material at
the level of the effective energy, it is convenient to identify it
with the ($\G$-)limit of a sequence of mixtures. Let $\{\chi_k\}$ be
a sequence in $\X(\O)$, we say that $W^*\in\F(\a,p,\O)$ is the
\emph{effective energy density} associated to
$\lbrace(W^{(i)},\chi^{(i)}_k)\rbrace_{i=1,\ldots,N}$ if the
functional $F:W^{1,p}(\O;\Rb^m)\rightarrow[0,+\infty)$ defined by
\begin{equation*}
F(u):=\int_\O W^*(x,\nabla u)\, dx
\end{equation*}
is the $\G$-limit of the sequence $\{F_{\chi_k}\}$. The sequence $
\{\chi_k\}$ is referred as the \emph{microstructure} (or
micro-geometry) of $W^*$. If the effective energy density does not
depend on $x$ we say that it is \emph{homogeneous} while if, for some
$\theta=(\theta^{(1)},\ldots,\theta^{(N)}) \in [0,1]^N$ with
$\sum_{i=1}^N\theta^{(i)}=1$, the equality
$$\med_\O \chi_k\,dx=\theta$$
holds for all $k\in\Nb$, we say that the micro-geometry $\{\chi_k\}$
has \emph{fixed volume fractions}
$\theta^{(1)},\ldots,\theta^{(N)}$.\\

The next lemma shows that every effective energy density can be
generated by a micro-geometry with fixed volume fractions.
\begin{lemma}\label{frazioni}
Let $\{\chi_k\}$ be a sequence in $\X(\O)$ and let
$W^*\in\F(\a,p,\O)$ be the effective energy density associated to
$\lbrace(W^{(i)},\chi^{(i)}_k)\rbrace_{i=1,\ldots,N}$. Suppose that
\begin{equation*}
\chi_k \xrightharpoonup[]{*} \theta \text{ in } L^\infty(\O;[0,1]^N)
\end{equation*}
and set $\overline{\theta}:=\med_\O \theta\,dx$. Then there exists
another sequence $\{\tilde{\chi}_k\}$ in $\X(\O)$ satisfying
\begin{equation*}
\tilde \chi_k \xrightharpoonup[]{*} \theta \text{ in }
L^\infty(\O;[0,1]^N), \quad \med_\O\tilde{\chi}_k\,
dx=\overline{\theta} \quad\text{for all } k\in\Nb\end{equation*} and
such that $W^*$ is the effective energy density associated to
$\lbrace(W^{(i)},\tilde{\chi}^{(i)}_k)\rbrace_{i=1,\ldots,N}$.
\end{lemma}
\begin{proof}For every $k \in \Nb$, let
$\theta_k:=(\theta^{(1)}_k,\ldots,\theta^{(N)}_k)$ where
$$\theta^{(i)}_k:=\med_\O \chi^{(i)}_k \, dx, \quad \text{for every }
i\in \{1,\ldots,N\}.$$ Let $I_k^1$, $I_k^2$ and $I_k^3$ be three
disjoint subsets of indexes in $\{1,\ldots,N\}$ such that
$$\left\{
\begin{array}{l}
\overline{\theta}^{(i)} > \theta^{(i)}_k \quad \text{if }i\in I_k^1,\\
\overline{\theta}^{(i)} = \theta^{(i)}_k \quad \text{if }i\in I_k^2,\\
\overline{\theta}^{(i)} < \theta^{(i)}_k \quad \text{if }i\in I_k^3.
\end{array}\right.$$
We denote $P_k^{(i)}:=\{\chi_k^{(i)}=1\}$ so that
$\theta^{(i)}_k=\LL^n(P_k^{(i)})/\LL^n(\O)$. For every $i \in
I_k^2$, we define $\tilde P^{(i)}_k:= P^{(i)}_k$ and $\tilde
\chi^{(i)}_k:=\chi_{\tilde P^{(i)}_k}$. Let $R>0$ be such that
$\O\subseteq B_R(0)$. Fix $i \in I_k^3$ and consider the function
$\ell_k(\rho):=\LL^n(B_\rho(0) \cap P^{(i)}_k)/\LL^n(\O)$; $\ell_k$
is a continuous and nondecreasing function satisfying $\ell_k(0)=0$
and $\ell_k(R)=\theta^{(i)}_k$. Hence one can find a radius $\rho_k
\in (0,R)$ such that $\ell_k(\rho_k)=\overline{\theta}^{(i)}$ and,
for any $i \in I_k^3$, we set $\tilde P^{(i)}_k:= P^{(i)}_k \cap
B_{\rho_k}(0)$ and $\tilde \chi^{(i)}_k:=\chi_{\tilde P^{(i)}_k}$.
It remains to treat the indexes $i \in I_k^1$. Set
$$Z_k:= \bigcup_{i\in I_k^3}
P_k^{(i)} \setminus \tilde P_k^{(i)};$$ Since
$$\sum_{i=1}^N\overline{\theta}^{(i)} = \sum_{i=1}^N\theta^{(i)}_k=1$$
it follows that
$$\LL^n(Z_k) =  \sum_{i\in I^3_k}\LL^n(\O) \big( \theta^{(i)}_k -
\overline{\theta}^{(i)} \big)=\sum_{i \in I_k^1} \LL^n(\O) \big(
\overline{\theta}^{(i)} - \theta_k^{(i)}\big).$$ Hence the measure
of the set $Z_k$ (that we have removed) is equal to the sum of the
measures of the sets that we need to add to each $P_k^{(i)}$ (for $i
\in I_k^1$) to get a larger set $\tilde P_k^{(i)}$ with Lebesgue
measure $\overline{\theta}^{(i)}\LL^n(\O)$. Thus we have enough room
to find disjoint measurable sets $\{Z_k^{(i)}\}_{i \in I_k^1}
\subseteq Z_k$ such that
$\LL^n(Z_k^{(i)})=\LL^n(\O)\big(\overline{\theta}^{(i)}-\theta^{(i)}_k\big)$
and it suffices to define $\tilde P_k^{(i)}:=P_k^{(i)} \cup
Z_k^{(i)}$ and $\tilde \chi^{(i)}_k:= \chi_{\tilde P_k^{(i)}}$ for
all $i \in I_k^1$. As
\begin{eqnarray}\label{misure}
\LL^n\left(\bigcup_{i=1}^N \big(P_k^{(i)} \triangle \tilde
P_k^{(i)}\big)\right) & = & \sum_{i \in I_k^1}\LL^n(\O) \big(
\overline{\theta}^{(i)} - \theta_k^{(i)}\big) + \sum_{i\in
I^3_k}\LL^n(\O) \big( \theta^{(i)}_k - \overline{\theta}^{(i)}
\big)\nonumber\\
& = & \LL^n(\O) \sum_{i=1}^N|\overline\theta^{(i)} - \theta_k^{(i)}|
\xrightarrow[k \to +\infty]{} 0,
\end{eqnarray}
one immediately gets that $\tilde \chi_k \xrightharpoonup[]{*}
\theta$ in $L^\infty(\O;[0,1]^N)$.

It remains to show that the $\G$-limit is unchanged upon replacing
$\chi_k$ by $\tilde \chi_k$. Extracting a subsequence (not
relabeled) if necessary, one may assume without loss of generality
(see Remark \ref{G} (ii-iii)) that there exists a function $\tilde
W\in\F(\a,p,\O)$ such that
$$\int_U \tilde W(x,\nabla u)\, dx=\G\text{-}\lim_{k \to +\infty} \int_U W_{\tilde \chi_k}(x,\nabla
u)\, dx,$$ for every $U\in\Op(\O)$ and every $u \in
W^{1,p}(U;\Rb^m)$.  Let us show that $\tilde W(x,\xi)=W^*(x,\xi)$
for all $\xi \in \Rb^{m \times n}$ and a.e. $x \in \O$. Fix
$\xi\in\Rb^{m\times n}$ and let $\{u_k\} \subset W^{1,p}(U;\Rb^m)$
be a sequence such that $u_k \rightharpoonup \xi \cdot$ in
$W^{1,p}(U;\Rb^m)$ and
$$\int_U \tilde W(x,\xi)\, dx=\lim_{k \to +\infty} \int_U W_{\tilde \chi_k}(x,\nabla
u_k)\, dx.$$ Thanks to the Decomposition Lemma (see
\cite[Lemma~1.2]{FMP}) one can assume without loss of generality
that the sequence $\{\nabla u_k\}$ is $p$-equi-integrable. By the
$p$-growth condition (\ref{pgrowth}) and (\ref{misure}), we get that
$$\int_U |W_{\tilde \chi_k}(x,\nabla u_k) - W_{\chi_k}(x,\nabla u_k)| \, dx \leq 2\a_3 \int_{U \cap
\bigcup_{i=1}^N\big( P_k^{(i)} \triangle \tilde P_k^{(i)}\big)}
(1+|\nabla u_k|^p)\, dx \to 0.$$ Hence
$$\int_U \tilde W(x,\xi)\, dx = \lim_{k \to +\infty} \int_U
W_{\chi_k}(x,\nabla u_k)\, dx \geq \int_U W^*(x,\xi)\, dx$$ and the
opposite inequality follows from a similar argument. Since the
equality
$$\int_U \tilde W(x,\xi)\, dx = \int_U W^*(x,\xi)\, dx$$
holds for every open subset $U$ of $\O$, we deduce that $\tilde
W(x,\xi)=W^*(x,\xi)$ for a.e. $x \in \O$.
\end{proof}

In order to state the localization result for mixtures, namely
Theorem \ref{localiz}, it will be more convenient to introduce the
set of functions $\Gt(\O)$ which is made of all possible homogeneous
energy densities obtained as $\G$-limits of mixtures of
$W^{(1)},\ldots,W^{(N)}$ and having asymptotically homogeneous
microstructure.

\begin{definition}
Let $\theta=(\theta^{(1)},\ldots,\theta^{(N)})\in[0,1]^N$ be such
that $\sum_{i=1}^N\theta^{(i)}=1$. We define $\Gt(\O)$ as the set of
all functions $W^*\in\F(\a,p)$ for which there exists a sequence
$\{\chi_k\}$ in $\X(\O)$ with the properties that
\begin{equation*}
\chi_k \xrightharpoonup[]{*} \theta \text{ in } L^\infty(\O;[0,1]^N)
\end{equation*}
and that $W^*$ is the effective energy density associated to
$\lbrace(W^{(i)},\chi^{(i)}_k)\rbrace_{i=1,\ldots,N}$.
\end{definition}
\begin{lemma}\label{dominio}
The set $\Gt(\O)$ is independent of the open set $\O\subseteq\Rb^n$.
\end{lemma}
\begin{proof}
Consider $\O$ and $\O'$ two arbitrary open subsets of $\Rb^n$ and
assume that $W^*\in \Gt(\O)$. Let us show that it belongs to
$\Gt(\O')$ too. From the definition of $\Gt(\O)$, there exists a
sequence of characteristic functions $\{\chi_k\}$ in $\X(\O)$ such
that $\chi_k \xrightharpoonup[]{*} \theta$ in $L^\infty(\O;[0,1]^N)$
and
$$\G\text{-}\lim_{k \to +\infty}\int_\O W_{\chi_k} (x,\nabla u)\,dx = \int_\O W^*(\nabla u)\,dx.$$
Let $a\in \O$ and $\rho>0$ be such that $a+\rho\, \O'\subseteq \O$.
It turns out that $\chi_k\xrightharpoonup[]{*} \theta$ in
$L^\infty(a+\rho\,\O';[0,1]^N)$ and, by a localization argument (see
Remark \ref{G} (iii)),
$$\G\text{-}\lim_{k \to +\infty}\int_{a+\rho\,\O'}W_{\chi_k} (x,\nabla u)\,dx=
\int_{a+\rho\,\O'}W^*(\nabla u)\,dx.$$ Define
$\tilde{\chi}_k(x):=\chi_k(a+\rho x)$ for all $x \in \O'$. Since
$\theta$ is constant, it is easily seen that
$\tilde\chi_k\xrightharpoonup[]{*} \theta$ in
$L^\infty(\O';[0,1]^N)$. Moreover a simple change of variables
together with the fact that $W^*$ does not depend on $x \in \O$
implies that
$$\G\text{-}\lim_{k \to +\infty}\int_{\O'} W_{\tilde{\chi}_k} (x,\nabla u)\,dx=
\int_{\O'} W^*(\nabla u)\,dx.$$
\end{proof}

Since $\Gt(\O)$ is independent of $\O$, we simply denote it from now
on by $\Gt$. We next show that every effective energy density can be
locally seen as homogeneous. Note that this result is independent of
any kind of convexity assumption on the densities
$W^{(1)},\ldots,W^{(N)}$.
\begin{theorem}\label{localiz}
Let $W^{(1)},\ldots, W^{(N)}\in\F(\a,p)$ and $\theta \in
L^\infty(\O;[0,1]^N)$ be such that $\sum_{i=1}^N\theta^{(i)}(x)=1$
a.e. in $\O$. Given $W^*\in\F(\a,p,\O)$, the following conditions
are equivalent:
\begin{itemize}
\item[i)] there exists a sequence $\{\chi_k\}$ in $\X(\O)$ such that
$\chi_k \xrightharpoonup[]{*} \theta$ in $L^\infty(\O;[0,1]^N)$ and
$W^*$ is the effective energy density associated to
$\lbrace(W^{(i)},\chi^{(i)}_k)\rbrace_{i=1,\ldots,N}$;

\vspace{5pt}
\item[ii)] $W^*(x,\cdot) \in \Gtx$ for a.e. $x \in \O$.
\end{itemize}
\end{theorem}

\begin{proof}
Thanks to Remark \ref{G} (iv) we can suppose without loss of
generality that
the functions $W^{(1)},\ldots, W^{(N)}$ are quasiconvex.\\
\textbf{(i)$\Rightarrow$(ii).} Let $\{\rho_j\} \searrow 0^+$ and fix
a point $x_0 \in \O \setminus Z$, where $Z \subseteq \O$ is the set
of Lebesgue measure zero given by Lemma \ref{rho} (with $f=W^*$),
which is also a Lebesgue point of $\theta$. Define
$\chi_{j,k}(y):=\chi_k(x_0 + \rho_j y )$, by Lemma \ref{rho} we have
that
\begin{eqnarray}\label{gagaconv}
&&\G\text{-}\lim_{j \to +\infty}\left( \G\text{-}\lim_{k \to
+\infty} \int_Q W_{\chi_{j,k}} (y,\nabla u(y))\, dy
\right)\nonumber\\
&&\hspace{3cm} = \G\text{-}\lim_{j \to +\infty}\left(
\G\text{-}\lim_{k \to +\infty} \int_Q W_{\chi_k} (x_0+\rho_j\, y,\nabla u(y))\, dy \right)\nonumber\\
&&\hspace{3cm} =  \G\text{-}\lim_{j \to +\infty}\int_Q W^*(x_0 + \rho_j y,\nabla u(y))\, dy\nonumber\\
&&\hspace{3cm} =  \int_Q W^*(x_0,\nabla u(y))\, dy.
\end{eqnarray}
Moreover, for every $\varphi \in C_c(Q)$ and any $i \in
\{1,\ldots,N\}$,
\begin{eqnarray}\label{lebpt}
&&\lim_{j \to +\infty}\lim_{k \to +\infty} \left|\int_Q
\Big(\chi^{(i)}_{j,k}(y) - \theta^{(i)}(x_0) \Big)\, \varphi(y)\, dy \right|\nonumber\\
&&\hspace{2cm} =  \lim_{j \to +\infty}\lim_{k \to +\infty} \left|
\med_{Q_{\rho_j}(x_0)} \left(\chi^{(i)}_k(y) -
\theta^{(i)}(x_0)\right)\,
\varphi\left(\frac{y-x_0}{\rho_j}\right)dy \right|\nonumber\\
&&\hspace{2cm} =  \lim_{j \to +\infty} \left| \med_{Q_{\rho_j}(x_0)}
(\theta^{(i)}(y) - \theta^{(i)}(x_0))\, \varphi\left(\frac{y-x_0}{\rho_j}\right)dy\right| \nonumber\\
&&\hspace{2cm} \leq \lim_{j \to +\infty}\|\varphi\|_{L^\infty(Q)}
\med_{Q_{\rho_j}(x_0)} |\theta^{(i)}(y) - \theta^{(i)}(x_0)|\, dy =0
\end{eqnarray}
since $x_0$ is a Lebesgue point of $\theta$. By density,
(\ref{lebpt}) remains true for any $\varphi \in L^1(Q)$. By
(\ref{gagaconv}), (\ref{lebpt}), the metrizable character of
$\G$-convergence for lower semicontinuous and coercive functionals
on $W^{1,p}(Q;\Rb^m)$  \cite[Corollary~10.22~(a)]{DM} and the
metrizability of the weak* convergence in $L^\infty(Q;[0,1])$, we
deduce through a standard diagonalization argument the existence of
a sequence $k_j \nearrow +\infty$ as $j \to +\infty$ (possibly
depending on $x_0$) such that $\tilde \chi_j:=\chi_{j,k_j}
\xrightharpoonup[]{*} \theta(x_0)$ in $L^\infty(Q;[0,1]^N)$ and
\begin{equation*}
\G\text{-}\lim_{j \to +\infty}
\int_Q W_{\tilde \chi_j} (y,\nabla u(y))\, dy = \int_Q
W^*(x_0,\nabla u(y))\, dy.
\end{equation*}

\vspace{5pt} \noindent \textbf{(ii)$\Rightarrow$(i).} By the Lusin
and the Scorza-Dragoni Theorems (see \cite{ET}), for every $k \in
\Nb$, there exists a compact set $K_k \subseteq \O$ satisfying
$\LL^n(\O \setminus K_k) < 1/k$ and such that
$\theta^{(1)},\ldots,\theta^{(N)}$ are continuous on $K_k$ and $W^*$
is continuous on $K_k \times \Rb^{m \times n}$. Let $h \in \Nb$,
split $\O$ into $h$ disjoint open sets $U_{r,h}$ such that
$\delta_h:=\max_{1\leq r \leq h}\text{diam}(U_{r,h}) \to 0$ as $h
\to +\infty$ and set
$$I_{h,k}:=\big\{ r \in \{1,\ldots,h\} : \LL^n(K_k \cap U_{r,h}) >0\big\}.$$ Hence for each
$r \in I_{h,k}$, one can find a point $x_{r,h}^k \in K_k \cap
U_{r,h}$ such that $W^*(x^k_{r,h},\cdot) \in \mathcal
G_{\theta(x^k_{r,h})}$. As a consequence of Lemma \ref{dominio},
there exist sequences
$\{\chi_{j}^{r,h,k}\}=\{((\chi_{j}^{r,h,k})^{(1)},\ldots,(\chi_{j}^{r,h,k})^{(N)})\}$
in $\X(U_{r,h})$ such that $\chi_{j}^{r,h,k} \xrightharpoonup[]{*}
\theta(x_{r,h}^k)$ in $L^\infty(U_{r,h};[0,1]^N)$ as $j \to +\infty$
and \begin{equation}\label{1127}\G\text{-}\lim_{j \to
+\infty}\int_{U_{r,h}} W_{\chi_{j}^{r,h,k}}(x, \nabla u)\,dx=
\int_{U_{r,h}} W^*(x_{r,h}^k,\nabla u)\,dx.\end{equation} Define now
\begin{equation*}\begin{split}
&W^*_{h,k}(x,\xi):=\sum_{r \in I_{h,k}} \chi_{U_{r,h}}(x)
W^*(x^k_{r,h},\xi)+ \chi_{U_{h,k}^c}(x) W^*(x,\xi)\\
\end{split}\end{equation*}
and
\begin{equation*}
\left\{\begin{array}{l}\ds (\chi_j^{h,k})^{(1)}(x):=\sum_{r\in
I_{h,k}} (\chi_j^{r,h,k})^{(1)}(x)\chi_{U_{r,h}}(x) +
\chi_{U^c_{h,k}}(x)\\
\ds (\chi_j^{h,k})^{(i)}(x):=\sum_{r\in I_{h,k}}
(\chi_j^{r,h,k})^{(i)}(x)\chi_{U_{r,h}}(x)\quad i \in
\{2,\ldots,N\},
\end{array}\right.\end{equation*}
where $U^c_{h,k}:= \O \setminus \bigcup_{r \in I_{h,k}} U_{r,h}$.
Then $\chi_j^{h,k}
:=((\chi_j^{h,k})^{(1)},\ldots,(\chi_j^{h,k})^{(N)}) \in \X(\O)$ and
for every $\varphi \in L^1(\O)$
\begin{eqnarray}\label{THETA}\lim_{j \to +\infty}
\int_\O (\chi_j^{h,k})^{(1)}(x) \, \varphi(x)\, dx & = & \lim_{j \to
+\infty} \sum_{r\in I_{h,k}} \int_{U_{r,h}}
(\chi_j^{r,h,k})^{(1)}(x) \, \varphi(x)\, dx
+\int_{U^c_{h,k}} \varphi(x)\, dx\nonumber\\
& = & \sum_{r \in I_{h,k}} \int_{U_{r,h}} \theta^{(1)}(x^k_{r,h}) \,
\varphi(x)\, dx+\int_{U^c_{h,k}} \varphi(x)\, dx,
\end{eqnarray}
while for each $i \in \{2,\ldots,N\}$
\begin{equation}\label{1050}\lim_{j \to +\infty}
\int_\O (\chi_j^{h,k})^{(i)}(x) \, \varphi(x)\, dx =\sum_{r \in
I_{h,k}} \int_{U_{r,h}} \theta^{(i)}(x^k_{r,h}) \, \varphi(x)\, dx.
\end{equation}
On the other hand, since $\|\theta^{(i)}\|_{L^\infty(\O)} \leq 1$ we
get that for any $i \in \{1,\ldots,N\}$
\begin{eqnarray}\label{THETA2}
&&\sum_{r\in I_{h,k}} \int_{U_{r,h}}
|\theta^{(i)}(x^k_{r,h})-\theta^{(i)}(x)| \, |\varphi(x)|
\,dx\nonumber\\
&&\hspace{3cm} \leq  \sum_{r\in I_{h,k}} \int_{U_{r,h}\cap K_k}
|\theta^{(i)}(x^k_{r,h})-\theta^{(i)}(x)|\, |\varphi(x)|  \, dx +
2 \int_{\O \setminus K_k}|\varphi|\, dx\nonumber\\
&&\hspace{3cm} \leq \o^{(i)}_k(\delta_h)\int_\O |\varphi|\, dx + 2
\int_{\O \setminus K_k}|\varphi|\, dx,
\end{eqnarray}
where $\o^{(i)}_k:[0,+\infty) \to [0,+\infty)$ is the modulus of
continuity of $\theta^{(i)}$ on $K_k$. Since $\theta^{(i)}$ is
uniformly continuous on $K_k$, it follows that $\o^{(i)}_k$ is a
continuous and nondecreasing function satisfying $\o^{(i)}_k(0)=0$.
Hence taking first the limit as $h \to +\infty$ and then as $k \to
+\infty$, we deduce from (\ref{THETA}), (\ref{1050}), (\ref{THETA2})
and the fact that $\LL^n(U_{h,k}^c) \leq \LL^n(\O \setminus K_k) \to
0$ (as $k \to +\infty$ uniformly with respect to $h \in \Nb$), that
\begin{equation}\label{THETA3}
\lim_{k \to +\infty}\lim_{h \to +\infty}\lim_{j \to +\infty} \int_\O
(\chi_j^{h,k})^{(i)}(x) \, \varphi(x)\, dx = \int_\O
\theta^{(i)}(x)\, \varphi(x)\, dx.
\end{equation}
Moreover, by virtue of the uniform continuity of $W^*$ on $K_k
\times B_R(0)$ for every $k\in \Nb$ and $R>0$, we get that
\begin{eqnarray*}
\int_\O \sup_{|\xi|\leq R} |W_{h,k}^*(x,\xi)-W^*(x,\xi)|\, dx & \leq
& \sum_{r \in I_{h,k}} \int_{U_{r,h} \cap K_k}
\sup_{|\xi|\leq R} |W^*(x^k_{r,h},\xi)-W^*(x,\xi)|\, dx\\
&& \hspace{1cm}+ 2\a_3(1+R^p) \LL^n(\O\setminus K_k)\\
& \leq & \o_{k,R}(\delta_h) \LL^n(\O)+ 2\a_3(1+R^p)
\LL^n(\O\setminus K_k)
\end{eqnarray*}
where $\o_{k,R}:[0,+\infty) \to [0,+\infty)$ is the modulus of
continuity of $W^*$ on $K_k \times B_R(0)$ which is a continuous and
nondecreasing function satisfying $\o_{k,R}(0)=0$. Hence taking the
limit as $h \to +\infty$ and then as $k\to +\infty$ we obtain that
for every $R>0$,
\begin{equation}\label{W*} \lim_{k \to +\infty}\lim_{h \to +\infty}
\int_\O \sup_{|\xi|\leq R} |W^*_{h,k}(x,\xi)-W^*(x,\xi)|\, dx =0.
\end{equation}
By (\ref{THETA3}), (\ref{W*}) and the fact that $L^1(\O)$ is
separable, one can apply a standard diagonalization technique to get
the existence of a sequence $h_k \nearrow +\infty$ as $k \to
+\infty$ such that, setting $\tilde \chi_j^k:=\chi_j^{h_k,k}$ and
$W^*_k:=W^*_{h_k,k}$, then for any $i \in \{1,\ldots,N\}$,
\begin{equation}\label{THETA4}
\lim_{k \to +\infty}\lim_{j \to +\infty} \int_\O (\tilde
\chi_j^k)^{(i)}(x) \, \varphi(x)\, dx = \int_\O \theta^{(i)}(x)\,
\varphi(x)\, dx.
\end{equation}
and for every $R>0$\begin{equation}\label{W*2} \lim_{k \to +\infty}
\int_\O \sup_{|\xi|\leq R} |W^*_k(x,\xi)-W^*(x,\xi)|\, dx =0.
\end{equation}
In particular, (\ref{W*2}) shows that (at least for a not relabeled
subsequence) for every $\xi \in \Rb^{m \times n}$, $W_k^*(\cdot,\xi)
\to W^*(\cdot,\xi)$ pointwise a.e. in $\O$ and by Remark \ref{G}
(v), we obtain that
$$\G\text{-}\lim_{k \to +\infty}\int_\O W_k^*(x,\nabla u)\,
dx=\int_\O W^*(x,\nabla u)\, dx.$$ Furthermore, thanks to
(\ref{1127}), one can show that
$$\G\text{-}\lim_{j \to +\infty}\int_\O W_{\tilde \chi_j^k}(x,\nabla u)\,dx=\int_\O
W_k^*(x,\nabla u)\, dx.$$ Indeed, the lower bound is immediate while
the construction of a recovery sequence can be completed using the
fact that, inside each $U_{r,h_k}$, there exists an optimal sequence
which matches $u$ on a neighborhood of $\partial U_{r,h_k}$ (see
Remark \ref{G} (vi)). In this way, we can glue continuously each
pieces on $U_{r,h_k}$ and extend by $u$ on the whole set $\O$ to
construct an optimal sequence. Thus we have that
\begin{equation}\label{1451}
\G\text{-}\lim_{k \to +\infty}\left(\G\text{-}\lim_{j \to +\infty}
\int_\O W_{\tilde \chi_j^k}(x,\nabla u)\,dx\right)=\int_\O
W^*(x,\nabla u)\, dx.
\end{equation}
Appealing once again to the metrizability of $\G$-convergence of
lower semicontinuous and coercive functionals on $W^{1,p}(\O;\Rb^m)$
and the metrizability of the weak* convergence in
$L^\infty(\O;[0,1])$, from (\ref{THETA4}) and (\ref{1451}) we
obtain, by a diagonalization process, the existence of a sequence
$j_k \nearrow +\infty$ as $k \to +\infty$ such that, upon setting
$\tilde \chi_k:=\tilde \chi_{j_k}^k$, then $\tilde \chi_k
\rightharpoonup \theta$ in $L^\infty(\O;[0,1]^N)$ and
$$\G\text{-}\lim_{k \to +\infty}\int_\O W_{\tilde \chi_k}(x,\nabla u)\,dx
= \int_\O W^*(x,\nabla u)\, dx,$$ which completes the proof of the
theorem. \end{proof}

\section{Characterization of zero level sets}\label{mesyoung}

\subsection{Effective energy densities}
\noindent The first aim of this section is to study the zero level
set of an effective energy density. This will be done by means of
gradient Young measures. As an application to the $G$-closure
problem in the convex case, one will show thanks to this result that
they may exist effective energy densities that cannot be exactly
reached by a periodic microstructure (Example \ref{tartar}).

\begin{theorem}\label{ponte}
Let $W^{(1)},\ldots,W^{(N)} \in \F(\a,p)$ and suppose that the
compact sets $\mathfrak{A}^{(i)}:=\lbrace \xi \in\Rb^{m\times n}
\;:\; W^{(i)}(\xi)=0\rbrace$ are pairwise disjoint. If $u\in
W^{1,p}(\O;\Rb^m)$ and $\theta \in L^\infty(\O;[0,1]^N)$ with
$\sum_{i=1}^N\theta^{(i)}(x)=1$ a.e. in $\O$, then the following
conditions are equivalent: \vspace{5pt}
\begin{itemize}
\item[i)] there exists a gradient Young measure $\nu\in\Y(\O;\Rb^{m\times n})$ such that
\begin{equation*}\left\{
\begin{array}{l}
\ds \text{supp }\nu_x \subseteq
\mathfrak{A}:=\bigcup_{i=1}^N\mathfrak{A}^{(i)} \quad\text{for
a.e. } x\in\O,\\
\ds \nu_x(\mathfrak{A}^{(i)})=\theta^{(i)}(x) \quad\text{for all }\;
i\in \{1,\ldots,N\} \quad\text{and a.e. }
x\in\O,\\[10pt]
\ds \overline{\nu}_x=\nabla u(x) \quad\text{for a.e. } x\in\O;
\end{array}\right.
\end{equation*}
\vspace{5pt}
\item[ii)] there exists a micro-geometry $\{\chi_k\}$ in $\X(\O)$ such that,
denoting by $W^*\in\F(\a,p,\O)$ the effective energy density associated
to $\lbrace(W^{(i)},\chi^{(i)}_k)\rbrace_{i=1,\ldots,N}$, then
\begin{equation*}\left\{\begin{array}{l}
\ds \med_\O\chi_k\,dx=\med_\O\theta\,dx \quad\text{for all } k\in\Nb,\\[5pt]
\ds \chi_k \xrightharpoonup[]{*} \theta \text{ in } L^\infty(\O;[0,1]^N),\\[8pt]
\ds W^*(x,\nabla u(x))=0 \quad\text{ for a.e. } x\in\O.
\end{array}\right.\end{equation*}
\end{itemize}
\end{theorem}

\vspace{1pt}
\begin{proof}
\textbf{(i)$\Rightarrow$(ii).} For sake of clarity, we divide the proof into four steps.\\

\noindent \textbf{Step 1}. By using Theorem \ref{fund young} (iii) and
the Decomposition Lemma \cite[Lemma~1.2]{FMP},
we can suppose that there exists a sequence $\{u_k\} \subset
W^{1,p}(\O;\Rb^m)$ such that $u_k \rightharpoonup u$ in
$W^{1,p}(\O;\Rb^m)$, $\{\nabla u_k\}$ is $p$-equi-integrable and
generates $\nu$ and $\mathrm{dist}(\nabla
u_k,\mathfrak{A})\rightarrow 0$ in measure. For any
$i\in\{1,\ldots,N\}$, we define the Borel subsets of $\Rb^{m \times
n}$
\begin{equation*}
\mathfrak{C}^{(i)}:=\lbrace \xi\in\Rb^{m\times n} \;:\;
\mathrm{dist}(\xi,\mathfrak{A}^{(i)})\leq\mathrm{dist}(\xi,\mathfrak{A}^{(h)})
\quad\text{for } h\neq i\rbrace
\end{equation*}
and, for every fixed $k\in\Nb$, the following measurable partitions
$\lbrace P^{(i)}_{k}\rbrace_{i=1,\ldots,N}$ of $\O$
\begin{equation*}\begin{split}
&P^{(1)}_{k}:=\lbrace x\in\O \;:\; \nabla u_k(x)\in\mathfrak{C}^{(1)}\rbrace,\\
&P^{(2)}_{k}:=\lbrace x\in\O \;:\; \nabla u_k(x)\in\mathfrak{C}^{(2)}\setminus\mathfrak{C}^{(1)}\rbrace,\\
&\quad\vdots\\
&P^{(N)}_{k}:=\lbrace x\in\O \;:\; \nabla
u_k(x)\in\mathfrak{C}^{(N)}\setminus\textstyle{\bigcup_{i=1}^{N-1}}
\mathfrak{C}^{(i)}\rbrace.
\end{split}\end{equation*}
Finally, we define the sequence $\{\chi_k\}$ in $\X(\O)$ by
$\chi_k:=(\chi_{P^{(1)}_{k}},\ldots, \chi_{P^{(N)}_{k}})$.
Extracting a subsequence (not relabeled) if necessary, we may assume
that there exists an effective energy density $W^*\in\F(\a,p,\O)$
associated to $\lbrace(W^{(i)},\chi^{(i)}_k)\rbrace_{i=1,\ldots,N}$ and
that the sequence $\{w_k\}:=\{(\chi_k,\nabla u_k)\}$
generates a Young measure $\mu\in\Y(\O;\Rb^N \times \Rb^{m \times n})$.\\

\noindent \textbf{Step 2}. Denoting by $\lbrace
e_1,\ldots,e_N\rbrace$ the canonical basis of $\Rb^N$, we claim that
for a.e. $x \in \O,$
\begin{equation}\label{somma}
{\mu_x}=\sum_{i=1}^N\delta_{e_i}\otimes ({\nu_x} \res\,
\mathfrak{A}^{(i)}).
\end{equation}
Given $x\in\O$ and $k\in\Nb$, we may find some $h\in\{1,\ldots,N\}$
so that $x\in P^{(h)}_{k}$ and thus
\begin{eqnarray*}
\mathrm{dist}\left(w_k(x),\bigcup_{i=1}^N(\{e_i\}\times\mathfrak{A}^{(i)})\right)
& \leq & \mathrm{dist}\bigl(w_k(x),\{e_h\}\times\mathfrak{A}^{(h)}\bigr)\\
&=&\mathrm{dist}\bigl(\nabla u_k(x),\mathfrak{A}^{(h)}\bigr)
=\mathrm{dist}\bigl(\nabla u_k(x),\mathfrak{A}\bigr).
\end{eqnarray*}
Then by Theorem \ref{fund young} (iii) we obtain that
\begin{equation*}
\mathrm{supp}\,\mu_x\subseteq\bigcup_{i=1}^N(\{e_i\}\times\mathfrak{A}^{(i)})
\quad\text{for a.e. }x\in\O.
\end{equation*}
Since the sets $\mathfrak{A}^{(1)},\dots,\mathfrak{A}^{(N)}$ are
pairwise disjoint, for any Borel set
$\mathfrak{B}\subseteq\mathfrak{A}^{(i)}$, we get that
\begin{equation*}
\mu_x(\{e_i\}\times\mathfrak{B})=\mu_x(\Rb^N\times\mathfrak{B})=\nu_x(\mathfrak{B})
\end{equation*}
and therefore
\begin{equation*}
{\mu_x}\res \,(\{e_i\} \times\mathfrak{A}^{(i)})=\delta_{e_i}\otimes
({\nu_x} \res\, \mathfrak{A}^{(i)}).
\end{equation*}
\\
\noindent \textbf{Step 3}. We prove that $W^*(x,\nabla u(x))=0$ for
a.e. $x\in\O$. By \eqref{somma} and applying Theorem~\ref{fund
young}~(ii) with $f(x,z,\xi)=\sum_{i=1}^N z^{(i)} W^{(i)}(\xi)$
(where $z^{(i)}$ is the $i$'th component of $z$), we have that
\begin{eqnarray*}
\int_\O W^*\bigl(x,\nabla u(x)\bigr)\, dx & \leq &
\lim_{k\rightarrow+\infty}\int_\O \sum_{i=1}^N\chi^{(i)}_k(x)W^{(i)}\bigl(\nabla u_k(x)\bigr)\, dx\\
&=&\int_\O\biggl(\int_{\Rb^N \times \Rb^{m \times n}}\sum_{i=1}^N
z^{(i)}W^{(i)}(\xi)\,d\mu_x(z,\xi)\biggr)\, dx\\
& = &
\int_\O\sum_{i=1}^N\biggl(\int_{\mathfrak{A^{(i)}}}W^{(i)}(\xi)\,d\nu_x(\xi)\biggr)\,
dx=0,
\end{eqnarray*}
where we used the fact that the sequence
$\{\nabla u_k\}$ is $p$-equi-integrable in the first equality.\\

\noindent \textbf{Step 4}. We prove that $\chi_k
\xrightharpoonup[]{*} \theta$ in $L^\infty(\O;[0,1]^N)$. Fix $i\in
\{1,\ldots,N\}$ and $\varphi\in L^1(\O)$, by \eqref{somma} and
Theorem \ref{fund young} (ii) (with $f(x,z,\xi)=\varphi(x)z^{(i)}$),
we have that
\begin{eqnarray*}
\lim_{k\rightarrow+\infty}\int_\O \varphi(x)\chi^{(i)}_k(x)\,dx & = &
\int_\O\biggl(\int_{\Rb^N \times \Rb^{m \times n}}\varphi(x)z^{(i)}d\mu_x(z,\xi)\biggr)dx\\
&  = &\int_\O
\varphi(x)\nu_x(\mathfrak{A}^{(i)})\,dx=\int_\O\varphi(x) \,
\theta^{(i)}(x)\, dx.
\end{eqnarray*}
Through Lemma \ref{frazioni} we can now modify the sequence
$\{\chi_k\}$ in order to get another one satisfying~(ii).

\vspace{6pt} \noindent \textbf{(ii)$\Rightarrow$(i).} Let $\{u_k\}$
be a recovery sequence such that $u_k \rightharpoonup u$ in
$W^{1,p}(\O;\Rb^m)$ and \begin{equation}\label{1043}\lim_{k \to
+\infty}\int_\O W_{\chi_k}(x,\nabla u_k)\, dx=\int_\O W^*(x,\nabla
u)\, dx=0.\end{equation} Up to a subsequence (not relabeled), we may
assume that the sequences $\{\nabla u_k\}$ and
$\{w_k\}:=\{(\chi_k,\nabla u_k)\}$ generate, respectively, the Young
measures $\nu\in\Y(\O;\Rb^{m\times n})$ and $\mu\in\Y(\O;\Rb^N
\times \Rb^{m\times n})$. According to Remark \ref{1049}, we get
that $\overline \nu_x=\nabla u(x)$ for a.e. $x \in \O$.

From (\ref{1043}) and Theorem \ref{fund young} (i), we get that
$$\int_\O\biggl(\int_{\Rb^N \times \Rb^{m \times n}}\sum_{i=1}^N
z^{(i)}W^{(i)}(\xi)\,d\mu_x(z,\xi)\biggr)\, dx \leq
\liminf_{k\rightarrow+\infty}\int_\O
\sum_{i=1}^N\chi^{(i)}_k(x)W^{(i)}\bigl(\nabla u_k(x)\bigr)\, dx =
0$$ and therefore by \cite[Lemma~3.3]{BJ2}
\begin{equation*}
\mathrm{supp}\,\mu_x\subseteq \bigcup_{i=1}^N
(\{e_i\}\times\mathfrak{A}^{(i)}) \quad\text{ for a.e. }x\in\O.
\end{equation*}
Thanks to the inequality
\begin{equation*}
\mathrm{dist}\left(w_k(x),\bigcup_{i=1}^N(\{e_i\}\times\mathfrak{A}^{(i)})\right)
\geq\mathrm{dist}\bigl(\nabla u_k(x),\mathfrak{A}\bigr),
\end{equation*}
we obtain from Theorem \ref{fund young} (iii) that
$\mathrm{supp}\,\nu_x\subseteq\mathfrak{A}$ for a.e. $x\in\O$.
Moreover, since the sets
$\mathfrak{A}^{(1)},\dots,\mathfrak{A}^{(N)}$ are pairwise disjoint,
as in steps 2 of the previous implication we have that
\begin{equation*}
\mu_x=\sum_{i=1}^N \delta_{e_i}\otimes (\nu_x \res\,
\mathfrak{A}^{(i)}).
\end{equation*}
Finally, for all $i\in \{1,\ldots,N\}$ and all $\varphi\in L^1(\O)$,
\begin{eqnarray*}
\int_\O\varphi(x)\nu_x(\mathfrak{A}^{(i)})\,dx &
= & \int_\O\biggl(\int_{\Rb^N \times \Rb^{m \times n}}\varphi(x)z^{(i)}d\mu_x(z,\xi)\biggr)dx\\
&=& \lim_{k\rightarrow+\infty}\int_\O \varphi(x)\,
\chi^{(i)}_k(x)\,dx =\int_\O \varphi(x)\, \theta^{(i)}(x)\,dx,
\end{eqnarray*}
and thus, by the arbitrariness of $\varphi$, it follows that
$\nu_x(\mathfrak{A}^{(i)})=\theta^{(i)}(x)$ for a.e. $x\in\O$.
\end{proof}

As a consequence of the localization result for effective energy
densities (Theorem \ref{localiz}) and of gradient Young measures
(Remark \ref{Y}), we deduce the following homogeneous version of
Theorem~\ref{ponte}.
\begin{corollary}\label{homog young}
Under the same hypothesis than that of Theorem \ref{ponte}, if
$A \in \Rb^{m \times n}$ and $\theta \in [0,1]^N$ with
$\sum_{i=1}^N\theta^{(i)}=1$, then the following conditions are
equivalent: \vspace{5pt}
\begin{itemize}
\item[i)] there exists a homogeneous gradient Young measure $\sigma\in\PP(\Rb^{m\times n})$ such that
\begin{equation*}
\text{supp }\sigma \subseteq \mathfrak{A},\quad \overline{\sigma}=A
\quad\text{and}\quad \sigma(\mathfrak{A}^{(i)})=\theta^{(i)}
\quad\text{for all }\; i\in \{1,\ldots,N\};
\end{equation*}
\item[ii)] there exists $W^*\in\Gt$ such that $W^*(A)=0$.
\end{itemize}
\end{corollary}

\vspace{9pt}
\subsection{Cell integrands}
\noindent We now characterize the zero level set of
$(W_{\chi})_{\cell}$ and point out its dependence on $\chi$ and on
the zero level sets of $W^{(1)},\ldots, W^{(N)}$. We refer to
\cite{Bar06} for a similar characterization of the zero level set of
$(W_{\chi})_{\hom}$. The result obtained has a lot of interesting
consequences because it will enable us to build new counter-examples
to the validity of the cell formula even in the quasiconvex case
(Example \ref{unico}) and to the continuity of the determinant with
respect to the two-scale convergence (Example \ref{no-unico}).

\begin{definition}
Given a measurable partition $\lbrace P^{(i)}\rbrace_{i=1,\ldots,N}$
of the unit cell $Q$ and a family of compact sets
$\lbrace\mathfrak{A}^{(i)}\rbrace_{i=1,\ldots,N}$, we define
$\mathfrak{A}_{\cell}$ as the set of matrices $\xi \in\Rb^{m\times
n}$ such that there exists $\phi\in W^{1,\infty}_{\per}(Q;\Rb^m)$
satisfying
\begin{equation*}
\xi+\nabla \phi(y)\in\mathfrak{A}^{(i)} \text{ for a.e. }\; y\in
P^{(i)} \quad \text{ and all }i\in\lbrace1,\ldots,N\rbrace.
\end{equation*}
We call $\mathfrak{A}_{\cell}$ the \emph{cell set} associated to
$\lbrace(\mathfrak{A}^{(i)},P^{(i)})\rbrace_{i=1,\ldots,N}$.
\end{definition}

\begin{lemma}\label{zero energy}
Let $\lbrace P^{(i)}\rbrace_{i=1,\ldots,N}$ be a measurable
partition of the unit cell $Q$ and define $\chi\in\X(Q)$ by
$\chi^{(i)}:=\chi_{P^{(i)}}$ for $i\in \{1,\ldots,N\}$. Assume
that the functions $W^{(1)},\ldots, W^{(N)}\neg\in\F(\a,p)$
are quasiconvex and also that the compact sets
$\mathfrak{A}^{(i)}:=\lbrace \xi \in\Rb^{m\times n} \;:\;
W^{(i)}(\xi)=0\rbrace$ are not empty. Then we have
\begin{equation*}
\mathfrak{A}_{\cell}=\Bigl\lbrace \xi\in\Rb^{m\times n} \;:\;
(W_\chi)_{\cell}(\xi)=0\Bigr\rbrace,
\end{equation*}
\vspace{3pt}
\noindent
\textit{i.e.}, the zero-level set of the cell
integrand $(W_\chi)_{\cell}$ associated to $W_\chi$ coincides with
the cell set $\mathfrak{A}_{\cell}$ associated to
$\lbrace(\mathfrak{A}^{(i)},P^{(i)})\rbrace_{i=1,\ldots,N}$.
\end{lemma}
\begin{proof}
We only prove the inclusion
$(W_\chi)_{\cell}^{-1}(0)\subseteq\mathfrak{A}_{\cell}$, since the
opposite one is immediate. If $(W_\chi)_{\cell}(\xi)=0$, then by
Remark \ref{minfcell} there exists $\phi \in
W^{1,p}_{\per}(Q;\Rb^m)$ such that
\begin{equation*}
\int_Q W_\chi\left(y, \xi+\nabla \phi(y)\right)dy=0
\end{equation*}
and so $W^{(i)}\left(\xi+\nabla \phi(y)\right)=0$ for a.e. $y\in
P^{(i)}$. As $\mathfrak{A}^{(1)},\ldots,\mathfrak{A}^{(N)}$ are
compact sets and $\xi + \nabla \phi(y) \in \bigcup_{i=1}^N
\mathfrak{A}^{(i)}$ for a.e. $y \in Q$, it follows that $\nabla \phi
\in L^\infty(\Rb^n;\Rb^{m \times n})$ and thus $\phi$ is Lipschitz
continuous. Consequently $\phi \in W^{1,\infty}_\per(Q;\Rb^m)$ which
proves that $\xi \in \mathfrak{A}_\cell$.
\end{proof}


\section{G-closure in the convex case}\label{gcl}

\noindent In this section we focus on effective energy densities
with periodic microstructure. Let
$\theta=(\theta^{(1)},\ldots,\theta^{(N)})\in[0,1]^N$ satisfying
$\sum_{i=1}^N\theta^{(i)}=1$, we define $P_\theta$ to be the set of
all functions $W^*\in\F(\a,p)$ for which there exists $\chi\in\X(Q)$
such that $\int_Q\chi\,dx=\theta$ and, defining the sequence
$\{\chi_k\}$ by
\begin{equation*}
\chi_k(y)=\chi\bigl(\langle k\,y\rangle\bigr) \quad\text{for } y\in Q,
\end{equation*}
$W^*$ is the effective energy density associated to
$\lbrace(W^{(i)},\chi^{(i)}_k)\rbrace_{i=1,\ldots,N}$. Since by the
Riemann-Lebesgue Lemma $\chi_k\xrightharpoonup[]{*} \theta$ in
$L^\infty(\O;[0,1]^N)$, we always have the inclusion
$P_\theta\subseteq\Gt$. In view of Theorem~\ref{BM},
\begin{equation*}\begin{split}
P_\theta =\Big\{W^* &\in \F(\a,p) \text{ :
 there exists }\chi \in \X(Q)\\
&\ds\text{such that } \int_Q \chi\, dy
=\theta \text{ and } W^*=(W_\chi)_\hom  \Big\}.
\end{split}\end{equation*}
Note that, in general, one cannot expect the set $P_\theta$ to be
already closed with respect to the pointwise convergence (see
Example \ref{tartar} below). Hence we also define
\begin{equation*}\begin{split}
\hspace*{10pt}G_\theta :=\Big\{W^* &\in \F(\a,p) \text{ :
there exists a sequence } \{W^*_k\} \,\text{ in }\, P_\theta\\
&\ds\text{such that }  W^*_k\rightarrow W^* \text{ pointwise}
\Big\}.
\end{split}\end{equation*}

Remark \ref{G} (v), the metrizability of $\G$-convergence of lower
semicontinuous and coercive functionals in $W^{1,p}(Q;\Rb^m)$ and a
standard diagonalization argument show that the sets $G_\theta$ and
$\Gt$ are closed for the pointwise convergence. Thus, since
$P_{\theta} \subseteq \Gt$, it follows that $G_\theta \subseteq
\Gt$. The following result states that, at least in the convex case,
this inclusion is actually an equality and thus it makes more
precise the result on locality of mixtures Theorem \ref{localiz}.
\begin{theorem}\label{BB}
Let $W^{(1)},\ldots, W^{(N)}\in\F(\a,p)$ be $N$ convex functions and
let $\theta \in L^\infty(\O;[0,1]^N)$ be such that
$\sum_{i=1}^N\theta^{(i)}(x)=1$ a.e. in $\O$. Given
$W^*\in\F(\a,p,\O)$, the following conditions are equivalent:
\vspace{5pt}
\begin{itemize}
\item[i)] there exists a sequence $\{\chi_k\}$ in $\X(\O)$ such that
$\chi_k \xrightharpoonup[]{*} \theta$ in $L^\infty(\O;[0,1]^N)$ and
$W^*$ is the effective energy associated to
$\lbrace(W^{(i)},\chi^{(i)}_k)\rbrace_{i=1,\ldots,N}$;

\vspace{5pt}
\item[ii)] $W^*(x,\cdot) \in G_{\theta(x)}$ for a.e. $x \in \O$.
\end{itemize}
\end{theorem}

\begin{proof}
From Theorem \ref{localiz}, it is sufficient to prove that
$G_{\theta}=\Gt$ for any fixed $\theta\in[0,1]^N$. Let $W^* \in
\Gt$. By Lemmas \ref{frazioni} and \ref{dominio}, there exists a
sequence $\{\chi_k\}$ in $\X(Q)$ such that $\int_Q\chi_k\,
dy=\theta$ for all $k\in\Nb$ and
\begin{equation}\label{hyp}
\G\text{-}\lim_{k \to +\infty} \int_Q W_{\chi_k} (y,\nabla u)\,dy = \int_{Q} W^*(\nabla u)\, dy.
\end{equation}
As a consequence of the classical property of convergence of
minimal values (see {\it e.g.} \cite{BD,DM})
\begin{eqnarray}\label{limsup}
W^*(\xi) & = & \min\left\{\int_Q W^*(\xi + \nabla \phi)\, dy :
\phi \in W^{1,p}_\per(Q;\Rb^m)\right\}\nonumber\\
 & = & \lim_{k \to +\infty} \min \left\{\int_Q W_{\chi_k}(y,\xi + \nabla \phi)\, dy :
\phi \in W^{1,p}_\per(Q;\Rb^m)\right\}\nonumber\\
 & = & \lim_{k \to +\infty} (W_{\chi_k})_\cell(\xi)
\end{eqnarray}
and the conclusion follows from the very definition of $G_\theta$
and the convexity assumption.\end{proof}

As the following example shows, even in the convex case the
inclusion $P_\theta \subseteq G_\theta$ is in general strict and
thus we cannot expect $P_\theta$ to be closed with respect to the
pointwise convergence.
\begin{example}\label{tartar}
Consider the following four matrices in $\Rb^{2\times2}$:
\begin{equation*}
A_1:=\dig(-1,-3), \;\;A_2:=\dig(-3,1), \;\;A_3:=\dig(1,3)
\;\;\text{and} \;\;A_4:=\dig\bigr(3,-1)
\end{equation*}
and let $\mathfrak{A}:=\{A_1,\ldots,A_4\}$. This set has a
peculiarity: despite the absence of rank-one connections, its
quasiconvex hull is not trivial. The relevance of this kind of sets
has been discovered independently by many authors and in different
contexts (\cite{AH,CT,MN,S,T2}, see also \cite{BFJK}). We refer to
\cite{Mul99bis} for a basic description and to \cite{KMS,Sz} for an
advanced analysis.
By \cite[Lemma~2.6~and~Example~d~in~Section~3.2]{Mul99bis}, if
$\theta:=(8/15,\,4/15,\,2/15,\,1/15)$, then
\begin{equation*}
\sigma:=\theta^{(1)}\delta_{A_1}+\theta^{(2)}\delta_{A_2}+\theta^{(3)}\delta_{A_3}+\theta^{(4)}\delta_{A_4}
\end{equation*}
is a homogeneous gradient Young measure. For $i\in\{1,\ldots,4\}$,
we consider the convex functions $W^{(i)}$ defined by
\begin{equation*}
W^{(i)}(\xi):=\snorm{A_i-\xi}^p.
\end{equation*}
Since $\overline{\sigma}=-I$, where $I:=\dig(1,1)$, from Corollary
\ref{homog young} and Theorem \ref{BB}, there exists some $W^*\in
G_\theta$ such that $W^*(-I)=0$. On the other hand, if $W^* \in
P_\theta$, then from Lemma \ref{zero energy}, there would exist
$\chi \in \X(Q)$ and $\phi \in W^{1,\infty}_\per(Q;\Rb^2)$ such that
$\nabla \phi(x)-I=A_i$ for a.e. $y \in P^{(i)}:= \{\chi^{(i)}=1\}$.
Defining $\varphi(y):=\phi(y)-y$, it follows that $\varphi$ is
Lipschitz continuous and that $\nabla \varphi(y)=A_i$ for a.e. $y
\in P^{(i)}$ which is impossible since the matrices $A_1,\ldots,A_4$
are not rank-one connected (see \cite{CK} or \cite{Kir}). Thus $W^*
\not\in P_\theta$.
\end{example}

\begin{remark}\label{Hausdorff}
Note that it is also possible to adapt the proofs of \cite{A} and
\cite{R} to our setting. Let us briefly explain how to proceed. We
first need to introduce a suitable metric space structure.
Specifically, let $\E_p$ be the space made of all continuous
functions $f : \Rb^{m \times n} \to \Rb$ such that there exists the
limit
$$\lim_{|\xi| \to +\infty}\frac{f(\xi)}{1+|\xi|^{p+1}}=0.$$
The space $\E_p$ endowed with the norm
$$\|f\|:=\sup_{\xi \in \Rb^{m \times n}}
\frac{|f(\xi)|}{1+|\xi|^{p+1}}$$ is a normed space. In fact, only
its metric structure will be useful for our purpose. Using the
Ascoli-Arzela Theorem and the fact that all the functions we are
considering here satisfy uniform $p$-growth and $p$-coercivity
conditions, one can show that the sets $G_\theta$ and $\Gt$ are
closed subsets of $\E_p$ and that $G_\theta$ is the closure of
$P_\theta$ in $\E_p$. Moreover, if $W^* \in \Gt$ and $\{\chi_k\}
\subset \X(Q)$ is a sequence of characteristic functions as in the
definition of $\Gt$, by convergence of minimizers, we always have
that $(W_{\chi_k})_\cell \to W^*$ pointwise and also in $\E_p$. The
idea now consists in using a Hausdorff convergence argument to
deduce directly from the previous property that $W^* \in G_\theta$.
We recall the definition of the Hausdorff distance between two
closed sets $G_1$ and $G_2$ in $\E_p$:
$$\dH(G_1,G_2):=\max\left\{ \sup_{g \in G_2}\inf_{f \in G_1}\|f-g\|\;  ,\;  \sup_{f \in G_1}\inf_{g \in
G_2}\|f-g\|\right\}.$$ To do that, we need to prove that $\theta
\mapsto G_\theta$ is continuous for the Hausdorff metric. Indeed, if
so, setting $\theta_k:=\int_Q \chi_k \, dx$, since
$(W_{\chi_k})_\cell \in P_{\theta_k} \subseteq G_{\theta_k}$ and
$\theta_k \to \theta$, it follows that $G_{\theta_k} \to G_\theta$
in the sense of Hausdorff. Hence remembering that
$(W_{\chi_k})_\cell \to W^*$ in $\E_p$, it implies that $W^* \in
G_\theta$.

To show that $G_\theta$ depends continuously on $\theta$, we prove
an estimate on $\dH(G_{\theta_1},G_{\theta_2})$ for any $\theta_1$
and $\theta_2 \in [0,1]^N$. This is done using a Meyers type
regularity result to the solutions of the minimization problem
(\ref{fcell}). Indeed, one can adapt the proof of {\it e.g.}
\cite[Theorem~3.1]{G} to get the existence of two universal
constants $c>0$ and $q>p$ (both depending only on $n$, $p$, $\a_1$,
$\a_2$ and $\a_3$) such that any solution $\varphi\in
W^{1,p}_{\per}(Q;\Rb^m)$ of (\ref{fcell}) has a higher integrability
property, namely that
$$\left( \int_Q |\nabla \varphi|^q\, dx \right)^{1/q} \leq
c(1+|\xi|).$$ Using this estimate together with the standard uniform
$p$-Lipschitz property satisfied by functions in $G_{\theta_1}$ and
$G_{\theta_2}$, one can show that
$$\dH(G_{\theta_1},G_{\theta_2}) \leq c' |\theta_1 -
\theta_2|^{(q-p)/q},$$ (where $c'>0$ depends only on $n$, $p$,
$\a_1$, $\a_2$ and $\a_3$) which proves the desired result.

This gives an alternative proof to the equality $\Gt=G_\theta$ and
stresses the fact that the Meyers type regularity result can be seen
as a stronger version of the Decomposition Lemma of \cite{FMP}.
\end{remark}


\section{Some counter-examples}\label{ce}

\noindent An example obtained in \cite{M} shows that the equality
between cell and homogenized integrands does not hold, even in the
quasiconvex case. Precisely, this result is obtained through a
rank-one laminated structure assembled by mixing two polyconvex
functions with $4$-growth in the case $n=m=2$. In this last section,
we present an alternative counter-example. Namely we show that when
$n=m=2$, there exist functions $W^{(1)}$ and $W^{(2)} \in \F(\a,p)$
(with $p>1$ arbitrary) and a measurable set $P \subseteq Q$ such
that, setting $W_\chi(y,\xi):=\chi_P(y)W^{(1)}(\xi)+\chi_{Q\setminus
P}(y)W^{(2)}(\xi)$, then $(W_\chi)_\hom < (W_\chi)_\cell$. We use an
argument based on the characterization of the zero level set stated
in Lemma~\ref{zero energy} to show that $W_{\cell}$ is not rank-one
convex.

\begin{example}\label{unico}
Consider the following matrices:
\begin{equation*}
O:=\dig(0,0), \;\;I:=\dig(1,1), \;\;A:=\dig(-1,1), \;\;B:=\dig(0,1)
\;\;\text{and} \;\;C:=\dig\bigr(0,1/2\bigl).
\end{equation*}
Let $W^{(1)}$ and $W^{(2)}\in\F(\a,p)$ be two quasiconvex functions
such that
\begin{equation}\label{zeri}
(W^{(1)})^{-1}(0)=\lbrace O,A\rbrace \quad\text{and}\quad
(W^{(2)})^{-1}(0)=\lbrace O,I\rbrace,
\end{equation}
and define the function
\begin{equation}\label{W_chi}
W_\chi(y,\xi):=\chi_P(y)W^{(1)}(\xi)+\chi_{Q\setminus
P}(y)W^{(2)}(\xi),
\end{equation}
where $P=(0,1/2)\times (0,1)$. Then, the function $(W_\chi)_\cell$
is not rank-one convex. Moreover, $W_\cell(C)>0$ while $W_\hom(C)=0$.
\end{example}

\begin{proof}
Obviously $(W_\chi)_\cell(O)=(W_\chi)_\hom(O)=0$. Let $\phi\in
W^{1,\infty}_\per(Q;\Rb^2)$ be defined by
$\phi(y):=(\snorm{y_1-1/2},0)$. It is immediate to check that
$\nabla\phi(y)=\chi_P(y)\dig(-1,0)+\chi_{Q\setminus P}(y)\dig(1,0)$ and
therefore $(W_\chi)_\cell (B)=(W_\chi)_\hom(B)=0$. Since
$(W_\chi)_\hom$ is rank-one convex (because it is quasiconvex),
necessarily $(W_\chi)_\hom(C)=0$. We shall prove that on the
contrary, $(W_\chi)_\cell(C)>0$. If not, by Lemma \ref{zero energy},
there exists $\phi\in W^{1,\infty}_\per(Q;\Rb^2)$ such that
\begin{equation*}
C+\nabla \phi(y)\in
\begin{cases}
\lbrace O,A\rbrace &\quad\text{for a.e.} \;y\in P;\\
&\\[-6pt]
\lbrace O,I\rbrace &\quad\text{for a.e.} \;y\in Q\setminus P.
\end{cases}
\end{equation*}
Since $\mathrm{rank}(I)=\mathrm{rank}(A)=2$, it follows from a classical
result on rigidity of Lipschitz functions (see {\it e.g.}
\cite[Proposition 1]{BJ}) that $\phi$ is an affine function on $P$
and $Q\setminus P$ and that either
$$C+\nabla \phi(y)=O$$
or
$$C+\nabla \phi(y)=\chi_P(y)A+\chi_{Q\setminus P}(y)I$$
for a.e. $y \in Q$. But in both cases we get a contradiction with
the fact that $\phi$ should be $Q$-periodic and thus
$(W_\chi)_\cell(C)\neq 0$.
\end{proof}

\begin{figure}
\begin{center}
\includegraphics[width=10cm]{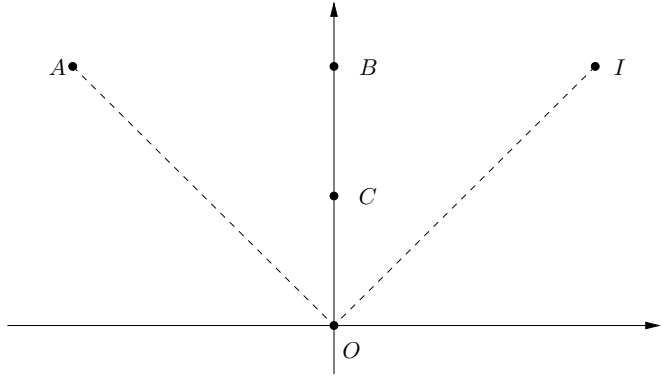}
\caption{A representation of the matrices $O$, $I$, $A$, $B$ and $C$ in $\Rb^2$,
identified with the set of the diagonal matrices.}
\end{center}
\end{figure}

\begin{remark}
Since the sets $\lbrace O,A\rbrace$ and $\lbrace O,I\rbrace$ are
quasiconvex, we can take as $W^{(1)}$ (resp. $W^{(2)}$) the
quasiconvexification of $D^{(1)}(\xi):=\mathrm{dist}^p(\xi,\lbrace
O,A\rbrace)$ (resp. $D^{(2)}(\xi):=\mathrm{dist}^p(\xi,\lbrace
O,I\rbrace)$).
\end{remark}

\begin{remark}
If a function $\phi\in W^{1,\infty}(Q;\Rb^2)$ is affine on $P$ and
\begin{equation*}
C+\nabla \phi(y)\in
\begin{cases}
\lbrace O,A\rbrace &\quad\text{for a.e.} \;y\in P;\\
&\\[-6pt]
\lbrace O,I\rbrace^\mathrm{co} &\quad\text{for a.e.} \;y\in
Q\setminus P,
\end{cases}
\end{equation*}
then it cannot be $Q$-periodic (the superscript $^\mathrm{co}$
denotes the convex hull of a set). By a straightforward modification
of Example \ref{unico}, we can take
$W^{(2)}(\xi):=\mathrm{dist}^p(\xi,\lbrace O,I\rbrace^\mathrm{co})$.
This proves that also by mixing a quasiconvex function and a convex
function, the equality $(W_\chi)_\cell=(W_\chi)_\hom$ could not
occur.
\end{remark}

\begin{remark}
If $p\geq2$, then we can take $W^{(1)}$ and $W^{(2)}$ polyconvex. We
recall that a function $W:\Rb^{2\times 2}\rightarrow\Rb$ is said to
be \emph{polyconvex} if there is a convex function
$V:\Rb^{2\times2}\times\Rb\rightarrow\Rb$ such that
$W(\xi)=V(\xi,\det(\xi))$ for all $\xi\in\Rb^{2\times 2}$. We refer
to \cite{D,Mul99bis} for more details.

\vspace{3pt} We define $V^{(1)}$ and
$V^{(2)}:\Rb^{2\times2}\times\Rb\rightarrow[0,+\infty)$ by
\begin{equation*}\begin{split}
&V^{(1)}(\xi,z):=\mathrm{max}\Bigl\lbrace\mathrm{dist}^p(\xi,\lbrace O,A\rbrace^\mathrm{co}),
\mathrm{dist}^{\frac{p}{2}}\bigl((\xi,z),\lbrace (O,0),(A,-1)\rbrace^\mathrm{co}\bigr)\Bigr\rbrace,\\
&V^{(2)}(\xi,z):=\mathrm{max}\Bigl\lbrace\mathrm{dist}^p(\xi,\lbrace
O,I\rbrace^\mathrm{co}),
\mathrm{dist}^{\frac{p}{2}}\bigl((\xi,z),\lbrace
(O,0),(I,1)\rbrace^\mathrm{co}\bigr)\Bigr\rbrace.
\end{split}\end{equation*}
Setting
\begin{equation}\label{Cex1}
W^{(1)}(\xi):=V^{(1)}(\xi,\mathrm{det}(\xi)) \quad\text{and}\quad
W^{(2)}(\xi):=V^{(2)}(\xi,\mathrm{det}(\xi)),
\end{equation}
we can verify that $W^{(1)}$ and $W^{(2)}\in\F(\a,p)$ for a suitable
$\a$ and that \eqref{zeri} holds. In particular $(W_\chi)_\hom(C)=0$
while $(W_\chi)_\cell(C)>0$.
\end{remark}

\vspace{7pt}

\noindent This last remark enables to underline the lack of
continuity of the determinant with respect to the two-scale
convergence. It is known (see {\it e.g.} \cite[Theorem~2.6 in
Chapter 4]{D}) that if $p>2$ and $\{u_k\}$ is a sequence in
$W^{1,p}(\O;\Rb^2)$, then
\begin{equation*}
u_k \rightharpoonup u \text{ in }W^{1,p}(\O;\Rb^2) \quad \text{
implies that }\quad \mathrm{det}(\nabla u_k) \rightharpoonup
\mathrm{det}(\nabla u) \text{ in }L^{p/2}(\O).
\end{equation*}
It would be interesting to ask whether this result still holds in
the two-scale convergence framework. The answer is negative as the
following example shows.

\begin{example}\label{no-unico}
Let $u(x):=Cx$, where $C=\text{diag}(0,1/2)$. By Theorem \ref{BM},
there exists a sequence $\{u_k\}$ in $W^{1,p}(\O;\Rb^2)$ such that
$u_k \rightharpoonup u$ in $W^{1,p}(\O;\Rb^2)$ and
\begin{equation*}
\lim_{k \to +\infty} \int_\O W_\chi\left(\left\langle \frac{x}{\e_k}
\right\rangle,\nabla u_k(x)\right)\, dx = \LL^n(\O)
(W_\chi)_\hom(C)=0,
\end{equation*}
where $W_\chi$ is defined by (\ref{W_chi}) with $W^{(1)}$ and
$W^{(2)}$ given by (\ref{Cex1}) for some $p>2$. By Remark
\ref{2scale} (i)-(iii), there exist a subsequence (not relabeled),
$v \in L^p(\O;W^{1,p}_\per(Q;\Rb^2))$ and $w \in L^{p/2}(\O \times
Q)$ such that $\nabla u_k \rightsquigarrow C + \nabla_y v$ and
$\mathrm{det}(\nabla u_k) \rightsquigarrow w$. Suppose now that
$w=\mathrm{det}(C + \nabla_y v)$, then by Remark \ref{2scale} (ii),
\begin{equation*}\begin{split}
\LL^n(\O)(W_\chi)_\cell(C)  \leq & \int_{\O \times Q}
W_\chi\big(y,C+\nabla_yv(x,y)\big)\, dx \,dy\\
 = &  \int_{\O \times Q}  \Bigl[\chi_P(y) V^{(1)}\big(C+\nabla_y
v(x,y),w(x,y)\big)\\
&  + \chi_{Q\setminus P}(y)
V^{(2)}\big(C+\nabla_y v(x,y),w(x,y)\big) \Bigr]\, dx\, dy\\
\leq &\liminf_{k \to +\infty}\int_\O
\Bigg[\chi_P\left(\left\langle\frac{x}{\e_k}\right\rangle\right)
V^{(1)}\big(\nabla u_k,\text{det}(\nabla u_k)\big)\\
&  + \chi_{Q \setminus P}
\left(\left\langle\frac{x}{\e_k}\right\rangle\right)
V^{(2)}\big(\nabla u_k,\text{det}(\nabla u_k)\big) \Bigg] dx\\
= & \lim_{k \to +\infty} \int_\O W_\chi\left(\left\langle
\frac{x}{\e_k} \right\rangle,\nabla u_k(x)\right)\, dx=0,
\end{split}
\end{equation*}
which is against the fact that $(W_\chi)_\cell(C)>0$.
Hence $w\neq\mathrm{det}(C+\nabla_y v)$.
This example can be connected to the recent investigation of the validity of the
div-curl lemma in two-scale convergence (see \cite{BCD,Vis}).
\end{example}

\bigskip

{\bf Acknowledgements. }The authors wish to thank Gianni Dal Maso
for having suggested this problem and for stimulating discussions on
the subject. They also gratefully acknowledge Gr\'egoire Allaire,
Andrea Braides, Adriana Garroni, Jan Kristensen and Enzo Nesi for
their useful comments and suggestions.

The research of J.-F. Babadjian has been supported by the Marie
Curie Research Training Network MRTN-CT-2004-505226 ``Multi-scale
modelling and characterisation for phase transformations in advanced
materials'' (MULTIMAT).

\end{document}